\documentclass[11pt]{amsart}
\usepackage{appendix}

\usepackage{mathrsfs}

\usepackage{a4wide}
\usepackage{amsmath,amssymb,amsfonts,amsthm,dutchcal}
\usepackage{tikz}
\usetikzlibrary{decorations.pathreplacing}
\usepackage{color}
\usepackage{graphicx}

\usepackage{cancel}

\usepackage[colorlinks=true,urlcolor=blue,
citecolor=red,linkcolor=blue,linktocpage,pdfpagelabels,
bookmarksnumbered,bookmarksopen]{hyperref}

\theoremstyle{plain}

\newtheorem{theorem}{Theorem}[section]
\newtheorem{lemma}[theorem]{Lemma}
\newtheorem{proposition}[theorem]{Proposition}

\newtheorem{remark}[theorem]{Remark}
\newtheorem{definition}[theorem]{Definition}
\theoremstyle{definition}
\numberwithin{equation}{section}

\def\R{\mathbb{R}}

\def\N{\mathbb{N}}
\def\H{\mathcal{H}}
\newcommand{\vphi}{\varphi} 
\newcommand{\la}{\langle}

\newcommand{\ra}{\rangle}
\renewcommand{\div}{\mathrm{div}}

\newcommand{\eps}{\varepsilon}

\newcommand{\pa}{\partial}

\def\beq{\begin{equation}}
\def\eeq{\end{equation}}

\title[Surface diffusion ]{Consistency for the Surface diffusion flat flow in three dimensions}

\author[M. Cicalese]
{M. Cicalese}
\address[Marco Cicalese]{
	Zentrum Mathematik - M7, Technische Universit\"at M\"unchen, Boltzmannstrasse 3, 85748 Garching, Germany	
}
\email[M. Cicalese]{cicalese@ma.tum.de}

\author[N. Fusco]
{N. Fusco}
\address[Nicola Fusco]{Dipartimento di Matematica e Applicazioni, Universit\`a di Napoli Federico II, Via Cintia, 80126 Napoli, Italy}
\email[N. Fusco]{n.fusco@unina.it}

\author[V. Julin]
{V. Julin}
\address[Vesa Julin]{
Jyv\"askyl\"an Yliopisto, Matematiikan ja Tilastotieteen Laitos, Jyv\"askyl\"a, Finland}
\email[V. Julin]{vesa.julin@jyu.fi
}

\author[A. Kubin]
{A. kubin}
\address[Andrea Kubin]{
	Jyv\"askyl\"an Yliopisto, Matematiikan ja Tilastotieteen Laitos, Jyv\"askyl\"a, Finland
}
\email[A. Kubin]{andrea.a.kubin@jyu.fi}

\begin{document}

\begin{abstract} 
We investigate the flat flow solution for the surface diffusion equation via the discrete minimizing movements scheme proposed by Cahn and Taylor \cite{CaTa}. 
We prove that in dimension three the scheme converges to the unique smooth solution of the equation, provided that the initial set is sufficiently regular. 
\end{abstract}

\maketitle

\section{Introduction}

In this paper we  study the evolution of smooth sets governed by the surface diffusion equation  
\begin{equation}
\label{eq:surf-diff}
V_t = \Delta_{\pa E_t} H_{E_t} \qquad \text{on } \, \pa E_t,
\end{equation}
where $V_t$ is the normal velocity, $ H_{E_t}$ the mean curvature and $ \Delta_{\pa E_t}$ the Laplace-Beltrami operator. We adopt the convention where the unit normal vector is outward and  $H_{E_t}$ is nonnegative for convex sets. It is immediate that this geometric flow preserves the volume of every component of the evolving sets, while decreasing the perimeter. 

The equation was proposed first by Mullins \cite{Mu} as a model to describe the development of surface grooves at the grain boundaries of a heated polycrystal. In the physically relevant three dimensional case it describes the evolution of interfaces in a variety of physical settings, including phase transition, epitaxial deposition and grain growth. The derivation of \eqref{eq:surf-diff} was later revisited by  Dav\`\i\, and Gurtin \cite{DaGu}  and by Cahn and Taylor \cite{CaTa}, who extended it to a more general mathematical context (see also \cite{GuJa}). In fact, 
Cahn and Taylor derived a more general equation of type 
\[
V_t = (\delta - \eps\Delta_{\pa E_t} )^{-1}  \Delta_{\pa E_t} H_{E_t} \qquad \text{on } \, \pa E_t.
\] 
Letting $\eps \to 0$  one obtains the surface-diffusion (with a constant $\delta$), while $\delta \to 0$  yields the volume preserving mean curvature flow 
\begin{equation}
\label{eq:VMCF}
V_t = - \eps  (H_{E_t}- \bar H_{E_t})  \qquad \text{on } \, \pa E_t.
\end{equation}
In \cite{CaTa} Cahn and Taylor also pointed out that the surface diffusion equation could be obtained as the sharp-interface limit ($\eps \to 0$) of the Cahn-Hilliard equation
\begin{equation}
\label{eq:Cahn-Hilliard}
\partial_t u = \div \big(m(u) \nabla( - \eps^2 \Delta u + W'(u) \big),
\end{equation}
where $W$ is a double well potential and $m(\cdot)$ is the mobility. This fact was formally proven in low dimensions $\R^n$, $n \leq 3$,  by Cahn, Elliot and Novik-Cohen in \cite{CaElNC}  with the mobility $m(u) = 1 - u^2$, the potential 
\[
W(u)  = \eps^\alpha \big((1+u) \log (1+u) +(1-u) \log (1-u) \big) +  (1-u^2),
\]
and by rescaling  the time $t \mapsto \eps^2 t$.  To the best of our knowledge the rigorous argument for this convergence is still missing, nor it is known that the limit yields a  weak solution to \eqref{eq:surf-diff}. We remark that by choosing the constant mobility $m(u) = 1$, the standard double well potential $W(u) = \tfrac14(1-u^2)^2$ and rescaling  the time $t \mapsto \eps t$  yields the Mullins-Sekerka flow as a limit as $\eps \to 0$, which was first formally derived by Pego  \cite{Pe}, and then  rigorously proven by  Alikakos, Bates and Chen \cite{AlBaCh}. We also mention the recent work  \cite{HoKiMuWa} where the authors consider a more general type Cahn-Hilliard equation and its sharp interface limit.  

For regular initial sets, the short time existence of a classical solution of  \eqref{eq:surf-diff}   was first proven in the planar case in \cite{BDR, EG, GigaIto} and later in every dimension by   Escher, Mayer and Simonett  \cite{EMS} for $C^{2,\alpha}$-regular initial sets. The result in \cite{EMS} holds even for immersed surfaces and the authors prove the global existence and exponential convergence for initial sets close to the ball. In the  flat torus $\mathbb{T}^n$ similar long time existence and convergence results near any stable critical set have been obtained in \cite{AFJM}, for $n=3$, and in \cite{DFM2023} for $n\geq4$. 

One expects that for generic initial sets the equation \eqref{eq:surf-diff}  forms singularities in finite time. Moreover, unlike the curve shortening flow, the flow under surface diffusion  for embedded curves may exhibit self-intersections in finite time as shown in \cite{GigaIto}. For immersed curves, singularities are also expected, as supported by the numerical simulation in \cite{EMS}.  We stress that in this paper, we require that the evolution in \eqref{eq:surf-diff}  is given by smooth and bounded sets and not immersed surfaces.

The construction of a weak solution to \eqref{eq:surf-diff}, which exists globally in time, presents a significant challenge due to the fourth-order nature of the equation and the absence of a maximum principle. A possible candidate for a weak solution, referred to as the flat curvature flow, was introduced by Cahn and Taylor in \cite{CaTa}. The starting point is the fact that  \eqref{eq:surf-diff} can  be seen as the $H^{-1}$-gradient flow of the area functional, see \cite{TaCa}. One may use the gradient flow structure to construct the solution using a minimizing movement scheme  similar to the one proposed by  Almgren, Taylor and Wang \cite{ATW} and  Luckhaus and St\"urzenhecker  \cite{LS} for the mean curvature flow. More precisely, given any initial bounded set of finite perimeter  $E_0 \subset \R^n$ and  a small time step $h>0$, one inductively defines $E_0^h= E_0$ and   $E_k^h$ for $k = 1,2,\dots$ as a minimizer of the functional
\begin{equation} \label{eq:def-scheme-intro}
P(F)+\frac{d_{H^{-1}}(F;E_{k-1}^h)^2}{2h},
\end{equation}
where  
\begin{equation} \label{eq:def-distance-intro}
d_{H^{-1}}(F;E)=\sup_{\|\nabla_{\pa E} f\|_{L^2(\pa E)}\leq1}\int_{\R^n}   f( \pi_{\pa E}(x))(\chi_F(x)-\chi_E(x))\,dx.
\end{equation}
Above, $\nabla_{\pa E}$ denotes the tangential gradient, $\chi_E$ the characteristic function of $E$ and $\pi_{\pa E}$ the projection on the boundary $\pa E$. Note that if $d_{H^{-1}}(F;E)$ is finite, then necessarily $|F| = |E|$, ensuring the scheme to be volume preserving. However, \eqref{eq:def-distance-intro} does not define a real distance,  nor does $d_{H^{-1}}(F;E) =0$ in general imply $F = E$.  
 Then Cahn and Taylor set $E^h(t) = E_k^h$ for $t \in [kh , (k+1)h)$ and define a flat flow as any cluster point of $E^h(t)$ as $h \to 0$. To the best of our knowledge, even the existence  of a minimizer of  the functional \eqref{eq:def-scheme-intro} is not known, let alone the convergence of the scheme as $h\to 0$ to a short time smooth solution of \eqref{eq:surf-diff} when the initial set $E_0$ is smooth. In order to highlight the difficulty of the problem,  let us consider the functional \eqref{eq:def-scheme-intro} in the flat torus $\mathbb{T}^n$ and the associated minimum problem
\[
\inf \bigg\{ P_{\mathbb{T}^n}(F)+\frac{1}{2h}\left(\sup_{\|\nabla_{\pa E} f\|_{L^2(\pa E)}\leq1}\int_{\mathbb{T}^n}   f( \pi_{\pa E}(x))(\chi_F(x)-\chi_E(x))\,dx \right)^2  \bigg\}.
\]
In the case of the strip  $E = (0,1)^{n-1} \times (0,\frac12) \subset \mathbb{T}^n$,  minimizers of the above problem are $E$ and $\mathbb{T}^n \setminus E$ for every $h$. This shows that at least in the case of flat torus, the scheme may fail to converge as $h \to 0$ unless the minimizer is chosen appropriately. 

The only result related to the existence of the flat flow  for the surface diffusion is contained in the thesis of Kin Yan Chung \cite{Chung}, where the author replaces the distance  \eqref{eq:def-distance-intro} with 
\begin{equation} \label{eq:def-distance-tilde}
\tilde d_{H^{-1}}(F;E)=\sup \{  \int_{\R^n}   f(x)(\chi_F(x)-\chi_E(x))\,dx : \,\, f : \R^n \to \R, \,\, \text{Lip}(f) \leq \Lambda, \,\,  \|\nabla_{\pa E} f\|_{L^2(\pa E)}\leq1\}. 
\end{equation}
 Using \eqref{eq:def-distance-tilde} in place of \eqref{eq:def-distance-intro} Chung shows the existence and continuity in time of the limiting flow of the minimizing movement scheme as $h \to 0$ and $\Lambda \to \infty$ (or by choosing $\Lambda = h^{-\alpha}$ for some $\alpha>0$). If $F$ and $E$ are smooth sets, whose Hausdorff distance is of order $h$, then it holds 
\[
\big|\tilde d_{H^{-1}}(F;E)  -  d_{H^{-1}}(F;E)  \big| \leq C h^2,
\]
for a constant that depends on $E$ and $\Lambda$, which means that \eqref{eq:def-distance-tilde} is formally a small perturbation of the original distance \eqref{eq:def-distance-intro}.  However, due to the double constraint it is not clear if  the  definition \eqref{eq:def-distance-tilde} in the minimizing movement scheme provides even formally  the surface diffusion \eqref{eq:surf-diff} as a limit as $h \to 0$.  Therefore the existence of a flat flow solution which is defined  globally in time and agrees with the classical solution for its existence time was still an open problem. 

The main result of this paper concerns the consistency of the scheme proposed by Cahn and Taylor in the physically relevant three dimensional case. More precisely, given a $C^5$-regular bounded initial set $E_0 \subset \R^3$, we fix a small $\delta >0$, set  $E_0^{h,\delta}= E_0$ and  define $E_k^{h,\delta}$ as a minimizer of  the following incremental minimum problem
\begin{equation}  \label{eq:def-scheme-intro-2}
\inf \Big{\{} P(F)+\frac{d_{H^{-1}}(F;E_{k-1}^{h,\delta})^2}{2h} : \,\, F \Delta E_{k-1}^{h,\delta} \subset \mathcal{N}_\delta(\pa  E_{k-1}^{h,\delta}) \Big{\}} ,
\end{equation}
where $\mathcal{N}_\delta(\Gamma)$ denotes the tubular neighborhood of a set $\Gamma \subset \R^3$ (see \eqref{def:tubular} for the definition). Due to the constraint condition $ F \Delta E_{k-1}^{h,\delta} \subset \mathcal{N}_\delta(\pa  E_{k-1}^{h,\delta})$ the existence of a minimizer of the above problem now easily follows from the direct methods of the Calculus of Variations. However, one of the key point of our proof is to show by means of quantitative geometric estimates, that any minimizer $ E_{k}^{h,\delta}$ of  \eqref{eq:def-scheme-intro-2} satisfies 
\begin{equation}  \label{eq:intro-3}
E_{k}^{h,\delta} \Delta E_{k-1}^{h,\delta} \subset \mathcal{N}_{\frac{\delta}{2}}(\pa  E_{k-1}^{h,\delta}),
\end{equation}
when $h$ is sufficiently small, provided $E_{k-1}^{h,\delta}$ is sufficiently regular. Therefore the additional constraint in \eqref{eq:def-scheme-intro-2} is merely a technical tool needed to  overcome the lack of coercivity  in the original functional \eqref{eq:def-scheme-intro}. We note that our main theorem holds also in the flat torus $\mathbb{T}^3$, in which case the constraint in \eqref{eq:def-scheme-intro-2} ensures that we choose a proper  minimizer of the problem.    We define the approximate flat flow $E^{h,\delta}(t)$ for all $t \geq 0$ setting $E^{h,\delta}(t) = E_k^{h,\delta}$ for $t \in [kh , (k+1)h)$ and define the constrained flat flow $(E^{\delta}(t))_{t \geq 0}$ as any cluster point of $E^{h,\delta}(t)$ as $h \to 0$ (see the precise definition in Section 4). 

Our main result is the short time regularity and the consistency of the minimizing movement scheme defined above. 
\begin{theorem}
\label{thm:thm2}
Assume $E_0  \subset \R^3$ is a bounded and $C^5$-regular set. There exist $\delta_0, \sigma_0$ and $C_0$ with the following property: for every $\delta < \delta_0$ there exists $h_0$ such that  the sets $E^{h,\delta}(t)$  are uniformly bounded in $H^4$, i.e., 
\[
\pa E^{h,\delta}(t) = \{ x + u^{h,\delta}(x,t) \nu_{E_0}(x) : \, x \in \pa E_0\} , \quad \|u^{h,\delta}(\cdot, t)\|_{H^4(\pa E_0)} \leq C_0
\]  
for all $t \in [0,\sigma_0]$ and $h \leq h_0$. Moreover, the functions $ u^{h,\delta}$ converge in $L^\infty(0,\sigma_0; C^{2,\alpha}(\pa E_0))$ for all $0<\alpha<1$ to a function $u^{\delta}$ such that the family  $(E^{\delta}(t))_{t \in [0,\sigma_0]}$ with  
\[
\pa E^{\delta}(t) = \{ x + u^{\delta}(x,t) \nu_{E_0}(x) : \, x \in \pa E_0\}
\]
is the unique smooth solution of \eqref{eq:surf-diff} starting from $E_0$. 
\end{theorem}


Note that, in particular,  the above result shows that the limit flow $(E^{\delta}(t))$ is the same for all $\delta < \delta_0$, which again indicates that the constraint in \eqref{eq:def-scheme-intro-2} is  irrelevant when $t \in [0,\sigma_0]$. We note that the novelty of Theorem \ref{thm:thm2} is the uniform higher regularity of the discrete flow, from where the consistency follows rather easily due to the uniqueness of the smooth solution.   In case the classical solution of  \eqref{eq:def-scheme-intro-2} is defined in a larger time interval $[0,T_0)$, $T_0 > \sigma_0$, we can prove the full consistency.   

\begin{theorem}
\label{thm:thm1}
Assume $E_0 \subset \R^3$ is as in Theorem \ref{thm:thm2} and assume the classical solution  $(E_t)$ of \eqref{eq:surf-diff} starting from $E_0$ exists in the time interval $[0,T_0)$. Then for every $T < T_0$ there exists $\delta(T)$ such that the  constrained flat flow   $(E^{\delta}(t))_{t \geq 0}$  agrees with $E_t$ in $[0,T]$ for all $\delta< \delta(T)$, i.e., $E_t =  E^{\delta}(t)$  for all  $t \in [0,T]$.  In addition, the sets $E^{h,\delta}(t)$ converge  to $E_t$  in $C^{2,\alpha}$ as $h \to 0$ uniformly in $[0,T]$, for all   $0<\alpha<1$.
\end{theorem}
Note that our theorems are stated in $\R^3$, but the results hold also in $\R^2$ following the same proof, with much less technical difficulties. 

As already mentioned, we are not aware of previous results on consistency of the flat flow solution or, in fact, any weak solution for the surface diffusion.  Construction of a  varifold solution via minimizing movement scheme to the Willmore flow, which is somewhat similar to the surface diffusion flow, is proven in \cite{BKS}, but even in this case  the consistency is not known.  The situation is better understood in the case of Mullins-Sekerka, where the flat flow solution constructed by   Luckhaus and St\"urzenhecker  \cite{LS} is known to be a weak solution of the associated equation in low dimension $n \leq 3$, due to results by R\"oger and Sch\"atzle \cite{Ro, Sch} and the consistency, or the weak-strong uniqueness,  is known to  hold in the planar case due to the recent result proven in  \cite{FHLS}. For the volume preserving mean curvature flow \eqref{eq:VMCF}, the third author with Niinikoski \cite{JN} proved the consistency of the flat flow solution constructed in \cite{MSS}. Finally we recall that the consistency of  the flat flow for the mean curvature flow follows from the comparison principle for the level set solutions \cite{Cha, CMP}. 


\subsection{Overview of the proof of Theorem \ref{thm:thm2}}

The crucial step in the proof of Theorem~\ref{thm:thm2} is to show that the speed of the discrete flow $E^{h,\delta}(t)$ is bounded in a suitable distance, i.e., the distance of the subsequent sets is linear in $h$. This is not too difficult in the case when the dissipation is given by an $L^2$-type distance as in the case of the mean curvature flow and volume preserving mean curvature flow \cite[Proposition 3.1]{JN}, where it follows from a direct comparison argument. This argument cannot be applied to  the case when the dissipation is given by the $H^{-1}$-distance, which contains less geometric information than the $L^2$-distance. In order to solve this issue we recall first that if a set $E$ satisfies the  uniform ball condition,  i.e., interior and exterior ball condition,  with radius $r$, then it is a $\Lambda$-minimizer of the perimeter in the sense that for all sets of finite perimeter $F$ it holds  (see \cite[Lemma 4.1]{AFM})
\[
P(E) \leq P(F) + \Lambda |E \Delta F|. 
\]
The above estimate is almost sharp in the sense that any  $\Lambda$-minimizer of the perimeter is $C^{1,\alpha}$-regular for any $\alpha \in (0,1)$ in $\R^n$ for $n \leq 7$. In Proposition \ref{prop:geometric} we prove that if  $E \subset \R^3$ satisfies the  uniform ball condition with radius $r$ and in addition $\|\Delta_{\pa E} H_{ E}\|_{L^2(\pa E)} \leq C$, then for all  sets of finite perimeter $F$ it holds
\[
P(E) \leq P(F) + \Lambda d_{H^{-1}}(F;E), 
\]
where $d_{H^{-1}}$ is defined in \eqref{eq:def-distance-intro}. This result can be seen as a higher order version of the $\Lambda$-minimality. Applying it to the minimization problem \eqref{eq:def-scheme-intro-2} we immediately obtain that the $d_{H^{-1}}$-distance between the sets $ E_{k-1}^{h,\delta}$ and $ E_{k}^{h,\delta}$ is linear in $h$ as long as they satisfy the regularity assumptions in  Proposition \ref{prop:geometric}.  

The second step is to show that we are able to keep the regularity assumed in  Proposition~\ref{prop:geometric} for a short time. To this aim we first prove in Proposition~\ref{prop:apirori-estimates} that the sets  $E_{k}^{h,\delta}$ are $(\Lambda,\alpha)$-minimizers of the perimeter (see Definition \ref{def:lambda-mini}) and use the well known  results of \cite{Tamanini1984} to deduce that when  the constraint parameter $\delta$  is small enough, the sets $E_{k}^{h,\delta}$ are uniformly $C^{1}$-regular with mean curvature bounded in $H^1$, and  can be written as a normal deformation over $E_0$ as described in Theorem \ref{thm:thm2}.  This, together with the linear  estimate $d_{H^{-1}}(E_{k}^{h,\delta}; E_{k-1}^{h,\delta})\leq Ch$, implies the estimate \eqref{eq:intro-3}, which is important since it means that the constraint in \eqref{eq:def-scheme-intro-2} is not active and that we may write the Euler-Lagrange equation.  We point out that  the argument for Proposition \ref{prop:geometric} can be generalized to higher dimension, but the proof of Proposition \ref{prop:apirori-estimates}  is constrained to low dimensions. The difficulties are similar to those of the Mullins-Sekerka flow \cite{Ro, Sch}, where the results for the flat flow hold in  $\R^3$. 

In order to get a uniform bound for  $\|\Delta_{\pa E_{k}^{h,\delta}} H_{ E_{k}^{h,\delta}}\|_{L^2}$, which we need for  Proposition~\ref{prop:geometric}, we use the fact that  the Euler-Lagrange equation \eqref{eq:Euler-Lagrange-2} is at each time step a fourth order elliptic equation. This allows us to  to upgrade the  $H^3$-regularity to $H^4$ (see Lemmas \ref{lem:upgrade-regularity} and \ref{lem:geometric-regularity}), with estimates which  however are not uniform in time.  The final step is to prove an iteration lemma (Lemma  \ref{lem:iteration}), where we use the parabolic nature of the flow and the estimates from the previous lemmas, to obtain  $H^4$-estimates uniform in time. With these high order estimates at hand, it is not surprising for the experts that  Theorems \ref{thm:thm2} and  \ref{thm:thm1} follow, but the proofs still require delicate technical arguments, which we chose to give in detail.

\section{Preliminaries}

\subsection{Notation}
We denote the Lebesgue measure of a set $E \subset \R^n$ by $|E|$.  We denote the distance function by $\text{dist}_E(x) := \inf_{y \in E}|x-y|$ and the signed distance function by $d_E:\R^{n} \rightarrow \R$, which is defined as
\[
d_E(x) := \begin{cases} \text{dist}_E(x) , \,\,&\text{for }\, x \in  \R^{n} \setminus E\\
- \text{dist}_{\R^n \setminus E}(x)  , \,\, &\text{for }\, x \in  E.  
\end{cases} 
\]
Clearly it holds $\text{dist}_{\pa E} = |d_E|$. If for a given  point $x \in \R^{n}$ there is a unique distance minimizer $y$ on 
$\pa E$ (that is $\text{dist}_{\pa E}(x) = |x-y|$), we call  $y$ the projection of $x$ onto $\pa E$ and set $y = \pi_{\pa E}(x)$.  Given $\sigma>0$ we define the tubular neighborhood of $\pa E$ as
\begin{equation}\label{def:tubular}
\mathcal{N}_\sigma (\pa E) := \{ x \in \R^n :  \text{dist}_{\pa E}(x) < \sigma \} .
\end{equation}
We denote the ball centred at $x$ with radius $r$ by $B_r(x)$ and by $B_r$ if it is centred at the origin.  We say that a set $E \subset \R^{n}$ satisfies the uniform ball condition (UBC) with a given radius $r >0$, 
if it simultaneously satisfies the exterior and interior ball condition with  radius $r$ at every boundary point. 
To be more precise,  for every 
$x \in \pa E$ there are balls $B_r(x_+)$ and $B_r(x_-)$  such that 
\[
B_r(x_+) \subset \R^{n}\setminus E, \quad B_r(x_-) \subset E \quad \text{and} \quad x \in \pa B_r(x_+) \cap \pa B_r(x_-).  
\]


We will mostly deal with regular and bounded sets $E \subset \R^{n}$.  We denote the outer unit normal by $\nu_E$. As usual, a bounded set $E \subset \R^{n}$ is said to be $C^{k,\alpha}$-regular, for some integer $k\geq1$ and an exponent $ 0<\alpha\leq1$, if its boundary is a $C^{k,\alpha}$-regular hypersurface, i.e., for every  $x \in \pa E$ we may find a radius $r>0$ and a function $f \in C^{k,\alpha}(\R^{n-1})$ such that, up to a rotation, 
\[
 E \cap B_r(x) =   \{(y',y_n) \in \R^n :  y_n < f(y')\}  \cap B_r(x).
\]
We say $E$ is smooth if $k = \infty$. If $E$ satisfies the UBC with radius $r_0$,  then by  \cite[Proposition~2.7] {JN} (see also \cite[Theorem 2.6]{Dal}  and \cite[Proposition~2.1]{MoMo})  the above holds  for $r = r_0/2$, i.e., up to a rotation
\begin{equation}\label{eq:UBC-graph-function}
 E \cap B_{r_0/2}(x) =   \{(y',y_n) \in \R^n :  y_n < f(y')\}  \cap B_{r_0/2}(x),
\end{equation}
and the function $f$ satisfies  $\|f\|_{C^{1,1}(B_{r_0/2}^{n-1}(x'))} \leq \frac{10}{r_0}$, assuming $r_0 \leq 1$. In particular, the radius is independent of the point $x$. Here by $B_{r}^{n-1}(x')$ we denote a ball in $\R^{n-1}$. We say that $E$ is uniformly $C^{k,\alpha}$-regular  for $k \geq 2$ with constants $r_0, C_0$, if it satisfies the UBC with radius $r_0$ and the function in  \eqref{eq:UBC-graph-function} satisfies $\|f\|_{C^{k,\alpha}(B_{r_0/2}^{n-1}(x'))} \leq C_0$. The notions of $C^k$-regular and uniformly $C^k$-regular function are defined similarly.

We define the matrix-valued function  $P_{\pa E} : \pa E \rightarrow \R^{n} \otimes \R^{n}$ by setting $P_{\pa E} = I - \nu_E\otimes
\nu_E$. For a given  point $x \in \pa E$ the map $P_{\pa E}(x)$ is the orthogonal projection onto the tangent hyperplane $T_x \pa E := \langle \nu_E (x) \rangle^\perp$. 

For a given map  $X \in C^1(\R^{n};\R^m)$ we define its tangential differential on $\pa E$ as
a matrix-valued function $\nabla_{\pa E} X : \pa E \rightarrow \R^m \otimes \R^{n}$ by setting
\[
\nabla_{\pa E} X :=\nabla X P_{\pa E} = \nabla X - (\nabla X   \nu_E) \otimes  \nu_E
\]
and, by slight abuse of notation,  the tangential gradient   of a real valued function $u : \R^n \to \R$ by $\nabla_{\pa E} u =P_{\pa E} \nabla u $. 
In the case $m=n$, the tangential divergence of $X$  is defined as $\div_{\pa E} X := \text{Tr}(\nabla_{\pa E} X)$.  In order to define the tangential Jacobian we denote the matrix associated with the natural inclusion by $I_x : T_x \pa E \to \R^n$, $I_x \tau = \tau$ for $\tau \in T_x \pa E$, and define the tangential Jacobian as $J_{\pa E} X(x):=  \sqrt{\det \big((\nabla_{\pa E} X(x) \, I_x)^T(\nabla_{\pa E} X (x)\, I_x)  \big)}$.  If $E$ is $C^2$-regular, we may define the tangential differential for maps  $X : \pa E \to \R^m$ by extending them in a tubular neighborhood of $\pa E$  as $\tilde X := X \circ \pi_{\pa E}$ and defining $\nabla_{\pa E}X := \nabla_{\pa E} \tilde X$ (see the discussion on the regularity of the projection below). 

If $E$ is $C^2$-regular, we define the symmetric matrix associated with the  second fundamental form as  $B_E := \nabla_{\pa E} \nu_E$. Note that this is not the standard  definition for the second fundamental form as $B_E$ is an $n\times n$-matrix with one eigenvalue always zero, but for our purposes this is the most convenient choice  of definition.   The mean curvature can be then written as  $H_E= \text{Tr}(B_E) = \div_{\pa E} \nu_E$ and the Gaussian curvature  $K_E(x) =\det \big(I_x^T B_E(x) I_x\big)$, where $I_x : T_x \pa E \to \R^n$ is the inclusion defined above.  We remark that for $E \subset \R^3$ we may write $B_E(x)$ for every $x \in \pa E$ as
\[
B_E(x) = \kappa_1 \, \tau_1 \otimes \tau_1 +  \kappa_2 \, \tau_2 \otimes \tau_2, 
\]
where $\kappa_1, \kappa_2$ are the principal curvatures and $\tau_1, \tau_2$ are the associated principal directions on the tangent hyperplane $T_x\partial E$. Therefore it is clear that our definition of $B_E$ is equivalent to the standard notion of the second fundamental form and, in particular, it holds $H_E = \kappa_1 +\kappa_2$ and $K_E = \kappa_1 \kappa_2$.  We use throughout the paper the fact that if $E$ satisfies the UBC with radius $r$, them $|B_E(x)| \leq \tfrac{n}{r}$.  
 We  have the divergence theorem  for all $X \in C^1(\pa E ;\R^{n})$ 
\[
\int_{\pa E} \div_{\pa E} X \, d \H^{n-1} = \int_{\pa E} H_E (X \cdot \nu_E) \, d \H^{n-1}.
\]
Finally, we define the Laplace-Beltrami  operator for a given $u \in C^2(\pa E)$ as 
\[
\Delta_{\pa E} u := \div_{\pa E}(\nabla_{\pa E} u).  
\]

If a bounded set $E$ satisfies the UBC with radius $r$, then the projection $\pi_{\pa E}$ is well defined in $ \mathcal{N}_{r}(\pa E)$ and can be written as 
\[
\pi_{\pa E}(x) = x - d_E(x) \nabla d_E(x). 
\] 
If in addition $E$ is $C^2$-regular, then $d_E \in C^2( \mathcal{N}_{r}(\pa E))$ and for any  $x \in \mathcal{N}_{r}(\pa E)$ its Hessian at $x$ is given by
\begin{equation}\label{eq:Hessian-distance}
\nabla^2d_E(x)=B_E(\pi_{\partial E}(x))(I+d_E(x)B_E(\pi_{\partial E}(x)))^{-1},
\end{equation}
see \cite[(2.33)]{JN}.
We note that   it is necessary to define the second fundamental form as the $n\times n$-matrix $
\nabla_{\pa E}\nu_E$ for the formula \eqref{eq:Hessian-distance} to hold. 
 Note also that the projection $\pi_{\pa E}$ is  $C^1$-regular and we may write its differential using \eqref{eq:Hessian-distance} and the fact that $\nabla d_E(x) = \nu_E(\pi_{\partial E}(x))$  
\begin{equation}\label{eq:derivative-proj}
\begin{split}
\nabla \pi_{\pa E}(x) &=  I - \nabla d_E \otimes \nabla d_E - d_E \nabla^2d_E   \\
			&= I - (\nu_E \circ \pi_{\pa E}) \otimes (\nu_E \circ \pi_{\pa E})  - d_E \, (B_E\circ \pi_{\partial E}) \big(I+d_E (B_E\circ \pi_{\partial E})\big)^{-1}. 
\end{split}
\end{equation}

We need also some notation from Riemannian geometry and as an introduction to the topic we refer to \cite{Lee}. Let us assume that $E \subset \R^{n}$ is a smooth and bounded set and set $\Sigma = \pa E$. Since $\Sigma$ is  embedded in $\R^{n}$  it has a natural metric $g$ induced by the Euclidean metric. Then $(\Sigma, g)$ is a Riemannian manifold and we denote the inner product of any two vectors  $X, Y$ in the tangent space $T_x \Sigma$ by $\la X, Y \ra$. Then, in local coordinates we have
\[
\la X,Y \ra = g(X,Y) = g_{ij} X^iY^j.
\]
We extend the inner product in a natural way for tensors, while we denote by $x \cdot y$ the inner product of two vectors in $\R^{n}$. We denote the space of smooth vector fields on $\Sigma$ by $\mathscr{T}(\Sigma)$ and by  a slight abuse of notation we denote the space of smooth $k$-th order tensor fields on $\Sigma$  by $\mathscr{T}^k(\Sigma)$.  We denote the Riemannian connection on  $\Sigma$ by $\bar \nabla$ and recall that  for a function $u \in C^\infty(\Sigma)$ the covariant derivative  $\bar \nabla u $ is a covector field defined for  $X  \in \mathscr{T}( \Sigma)$  as
\[
\bar \nabla u(X)  := \bar  \nabla_X u = X u,
\]
i.e., the derivative  of $u$ in the direction of $X$.  We denote the  covariant derivative  of  a smooth $k$-tensor field $F \in \mathscr{T}^k( \Sigma)$  by  $\bar \nabla F$, which is  a $(k+1)$-tensor field, for the definition see \cite{Lee},   and  the $k$-th order covariant derivative of a function $u$ on $\Sigma$ by $\bar \nabla^k u \in \mathscr{T}^k( \Sigma)$. 

We will always identify the tangent hyperplane $T_x \pa E := \langle \nu_E (x) \rangle^\perp$ and the tangent space $T_x \Sigma$ as defined in Differential Geometry.  With this convention  the tangential gradient of $u : \Sigma = \pa E \to \R$ is equivalent to its covariant derivative in the sense that for every tangent vector $\tau \in T_x \pa E$  it holds 
\[
\bar \nabla u(\tau)  = \nabla_{\pa E} u \cdot \tau. 
\]

Since $\Sigma$ is embedded in $\R^n$ we may define the second fundamental form and denote it by $B_\Sigma$. Note that in this setting $B_\Sigma(x)$  for $x \in \Sigma$ is understood as a bilinear form in $T_x\Sigma \times T_x\Sigma$.  

We define  the Sobolev space $W^{l,p}(\Sigma)$ in  a standard way for $p \in [1,\infty]$, see e.g. \cite{AubinBook2}, denote the Hilbert space $H^l(\Sigma) = W^{l,2}(\Sigma)$ and  define the associated norm for $f \in W^{l,p}(\Sigma)$ as
\[
\| f\|_{W^{l,p}(\Sigma)}^p = \sum_{k = 0}^l \int_\Sigma |\bar \nabla^k f|^p\, d \H^{n-1}
\]
and for $p = \infty$
\[
\| f\|_{W^{l,\infty}(\Sigma)} = \sum_{k = 0}^l \text{ess}\sup_{\!\!\!\!\!\!\!\!\!\!\!x \in \Sigma} |\bar \nabla^k f|.
\]
The above definition extends naturally to tensor fields.  We adopt the convention that $\| u\|_{H^0(\Sigma) }  = \| u\|_{L^2(\Sigma)}$ and  denote $\| u\|_{C^{m}(\Sigma)} = \| u\|_{W^{m,\infty}(\Sigma)}$ for $m \geq 0$. We remark that we may define the  space $W^{k,p}(\Sigma)$ for $k \geq 2$ as above assuming that $\Sigma$ (i.e., the set $E$ for which $\Sigma = \pa E$)  is $C^{k-1,1}$-regular.  We remark that we may define the first order Sobolev space for real valued functions by using the tangential derivatives and we will do this throughout the paper. However, when we need higher order Sobolev spaces or Sobolev space for tensor fields we need to rely on covariant derivatives.

Finally, we define the H\"older norm of a continuous function $u : \Sigma \to \R$ by 
\[
\| u\|_{C^\alpha(\Sigma)} =  \sup_{\substack{x \neq y \\ x,y \in \Sigma}} \frac{|u(y) - u(x)|}{d(y,x)^\alpha} + \|u\|_{L^\infty(\Sigma)},
\]
where $d(y,x)$ denotes the geodesic distance between $x$ and $y$ on $\Sigma$.  
We define the H\"older norm of a  tensor field $T \in \mathscr{T}^k(\Sigma)$ as in \cite{JL}
\[
\| T\|_{C^\alpha(\Sigma)} = \sup \{ \| T(X_1, \dots, X_k) \|_{C^\alpha(\Sigma)} : X_i \in \mathscr{T}(\Sigma) \, \text{ with } \, \|X_i \|_{C^1(\Sigma)} \leq 1\}.
\]


\subsection{Preliminary Results}

We begin by recalling the interpolation inequalities with Sobolev-norms on embedded surfaces. We use the result from \cite[Proposition 6.5]{Mantegazza2002} (see also \cite[Proposition 4.3]{DFM2023}) which states that under a curvature bound the standard interpolation inequalities hold with a uniform constant.   
\begin{proposition}
\label{prop:interpolation}
Assume $\|B_{\Sigma}\|_{L^{\infty}}, \H^{n-1}(\Sigma)\leq C_0$ and $\Sigma$ is $C^{m}$-regular for $m \geq 2$.  Then for integers $0\leq k < l \leq m$ and numbers $p \in [1,\infty)$,  $q,r \in [1,\infty]$   there is $\theta \in [k/l,1]$ such that for every $C^l$-regular  tensor field $T$ on $\Sigma$ it holds
\[
\|\bar \nabla^k T\|_{L^p(\Sigma)} \leq C \| T\|_{W^{l,q}(\Sigma)}^\theta \| T\|_{L^{r}(\Sigma)}^{1-\theta}
\]
for a constant $C=C(k,l,n,p,q,r,\theta,C_0)$,  provided that the following  condition is satisfied
\[
\frac{1}{p} = \frac{k}{n-1} + \theta \left( \frac{1}{q} - \frac{l}{n-1} \right) + \frac{1}{r}(1 - \theta).
\]
Moreover, if $f: \Sigma \to \R$ is a smooth function with $\int_{\Sigma_i} f\,d\H^{n-1}=0$ on every component $\Sigma_i$ of $\Sigma$, the above inequality can be written as 
\[
\|\bar \nabla^k f\|_{L^p(\Sigma)} \leq C \| \bar \nabla^l f\|_{L^{q}(\Sigma)}^\theta \| f\|_{L^{r}(\Sigma)}^{1-\theta}.
\]
\end{proposition} 
With the above assumptions on $\Sigma$, we also have the standard Sobolev embedding $\|f\|_{C^{\alpha}(\Sigma)} \leq C_p \|f\|_{W^{1,p}(\Sigma)}$ for $\alpha = 1- \frac{n-1}{p}$ and $p> n-1$. 
If $E \subset B_{R_0} \subset  \R^n$ is bounded and uniformly $C^{m,\beta}$-regular, for $m \geq 2$ with constants $r_0,C_0$, then we have also the classical  interpolation in H\"older norms, i.e.,  for $0 < \alpha < \beta \leq1$ and $0 \leq k \leq l \leq m$ it holds
\begin{equation} \label{eq:Holder-interpolate}
\|f\|_{C^{k,\alpha}(\Sigma)} \leq C \|f\|_{C^{l,\beta}(\Sigma)}^\theta\|f\|_{C^{0}(\Sigma)}^{1-\theta}, \qquad \text{for }\, \theta = \frac{k+\alpha}{l+\beta}, 
\end{equation}
where the constant depends on $r_0, R_0, C_0, k,l,\alpha,\beta $.  The above formula follows from the Euclidean result, see for instance \cite[Example 1.9]{Lunardi2018}.

The interpolation inequality in Proposition \ref{prop:interpolation} implies the following useful estimate.  The proof is standard and we refer to \cite[Proposition 2.3]{JN}, and note that the proof is similar as in the Euclidean case \cite[Proposition 3.7]{Taylor}. We denote the norm of an index vector  $\alpha \in \N^l$ by
\[
 |\alpha| := \alpha_1+\cdots + \alpha_l 
\]
\begin{lemma}
\label{lem:Leibniz}
Assume $\|B_{\Sigma}\|_{L^{\infty}}, \H^{n-1}(\Sigma)\leq C_0$ and $\Sigma$ is $C^{m}$-regular for $m \geq 2$, and assume $T_1, \dots, T_l$ are $C^m$-regular tensor fields.  Then for an index vector $\alpha \in \N^l$ with norm $|\alpha|  \leq  k \leq m$ it holds 
\[
\| |\bar \nabla^{\alpha_1} T_1| \cdots  |\bar \nabla^{\alpha_l} T_l| \|_{L^2(\Sigma)} \leq C \sum_{\sigma \in S_l} \|T_{\sigma(1)}\|_{L^\infty(\Sigma)} \cdots  \|T_{\sigma(l-1)}\|_{L^\infty(\Sigma)}  \,    \|T_{\sigma(l)}\|_{H^k(\Sigma)}, 
\]
where the sum is over the permutation group $S_l$.  In particular, 
\[
\| |\bar \nabla^{k} \la T_1 , T_2\ra  | \|_{L^2(\Sigma)} \leq C \|T_1\|_{L^\infty(\Sigma)}   \|T_2\|_{H^k(\Sigma)}    + C \|T_2\|_{L^\infty(\Sigma)}   \|T_1\|_{H^k(\Sigma)}  . 
\] 
\end{lemma}
We also need to bound the whole $\|u\|_{H^{2k}}$-norm with the norm $\| \Delta^k u\|_{L^{2}}$. Similarly we need to bound   $\|B_\Sigma\|_{H^{2k}}$-norm with the norm $\| \Delta^k H_{\pa E}\|_{L^{2}}$ This is given by the following lemma whose proof can be found e.g. in \cite[Lemma 2.5 \& Proposition 2.6]{JN}.
\begin{lemma}
\label{lem:hilbert-norm}
Assume $\Sigma = \pa E \subset \R^{n}$ is $C^{2k+1}$-regular with $ \H^{n-1}(\Sigma), \|B_\Sigma\|_{L^\infty} \leq C$. Then for all $f \in C^{2k}(\Sigma)$ it holds
\[
\begin{split}
&\|f\|_{H^{2k}(\Sigma)} \leq C_k\big( \|\Delta^k f\|_{L^{2}(\pa E)} +  (1+  \| \nabla_{\pa E} \Delta^{k-1}  H_{\pa E} \|_{L^{2}(\pa E)})  \|f\|_{L^{\infty}(\pa E)}\big),\\
&\|f\|_{H^{2k-1}(\Sigma)} \leq C_k\big( \| \nabla_{\pa E} \Delta^{k-1} f\|_{L^{2}(\pa E)} +  (1+ \| \Delta^{k-1} H_{\pa E} \|_{L^{2}(\pa E)}  )  \|f\|_{L^{\infty}(\pa E)}\big).
\end{split}
\]

Moreover, it holds 
\[
\begin{split}
&\|B_\Sigma \|_{H^{2k-2}(\Sigma)} \leq C_k (1 + \| \Delta^{k-1} H_{\pa E} \|_{L^{2}(\pa E)} ), \\
&\|B_\Sigma \|_{H^{2k-1}(\Sigma)} \leq C_k (1 + \| \nabla_{\pa E} \Delta^{k-1}  H_{\pa E} \|_{L^{2}(\pa E)} ).
\end{split}
\]
Here the constant $C_k$ depends only on $n$, $k$ and $C$.
\end{lemma}

The estimates for $\|f\|_{H^{2k}(\Sigma)}$ in Lemma \ref{lem:hilbert-norm} are sharp with respect to the norm on the curvature. However, since we assume that the second fundamental form is point-wise uniformly bounded the following holds under the assumptions of   Lemma \ref{lem:hilbert-norm}
\begin{equation}\label{eq:hilbert-norm}
\|f\|_{H^{2}(\Sigma)} \leq C\big( \|\Delta f\|_{L^{2}(\pa E)} + \|  f\|_{L^{2}(\pa E)} \big). 
\end{equation}
This follows from a similar argument as the Lichnerowicz's theorem \cite[Theorem 4.19]{AubinBook2} (see e.g. \cite[Remark 2.4 \& 2.5]{FJM3D}).

The above results hold for standard Sobolev-norms. In the definition of the minimizing movement scheme it appears a geometric $H^{-1}$-distance, which is defined for sets of finite perimeter $E, F$ as
\begin{equation}\label{def:akka-meno-uno}
d_{H^{-1}}(F;E)=\sup_{\|\nabla_{\pa E} f\|_{L^2(\pa E)}\leq1}\int_{\R^n}   f( \pi_{\pa E}(x))(\chi_F(x)-\chi_E(x))\,dx.
\end{equation}
Note that the projection $ \pi_{\pa E}$ is well-defined almost everywhere in $\R^n$. Indeed, this follows from the fact that the projection is unique at points where   the signed distance function $d_E$ is differentiable.  Note  that \eqref{def:akka-meno-uno} is not symmetric as it is defined using the projection on $E$. 

 It is  immediate that if  $d_{H^{-1}}(F;E)$ is finite then necessarily $|F| = |E|$. In general however, it is not clear what is the sufficient condition for the set $F$ to have finite distance in  \eqref{def:akka-meno-uno}. The situation becomes more clear if  we assume that $E$ satisfies the UBC with radius $r_0$. Indeed, then  the signed distance function $d_E$ is differentiable in $\mathcal N_{r_0}(\pa E)$, the projection is well defined there and we may define a map $\Phi_\tau : \pa E  \to \{ x \in \R^n : d_E(x) = \tau\}$, with $\tau \in (-r_0,r_0)$, as 
\begin{equation}
\label{def:mapping-Phi}
\Phi_\tau(x) = x + \tau \nu_{E}(x)
\end{equation}
which is a bijection. Moreover, the boundary $\pa E$ has finite number of components $\Sigma_1, ..., \Sigma_N$.  In this case we may characterize the sets  $F$ with finite $d_{H^{-1}}$-distance from $E$ when we constraint them to satisfy $F \Delta E \subset \mathcal N_{r_0}(\pa E)$. We state  this  as a lemma. We state the lemma in three dimensions in order to simplify the technicalities, but the result easily generalizes to any dimension.   

\begin{lemma}
\label{lem:distance-finite}
Assume that  $E \subset \R^3$ is $C^{2}$-regular,  $E \subset B_{R_0}$, satisfies  the UBC with radius $r_0$ and denote the components of $\pa E$ by $\Sigma_i$.  Let  $F$ be  a measurable set such that
\[
F \Delta E \subset \mathcal N_{r_0}(\pa E)
\]  
and define the function $\xi_{F,E} : \pa E \to \R$ as
\begin{equation}
\label{def:xi-original}
\xi_{F,E}(x) := \int_{-r_0}^{r_0}  (\chi_F( x + \tau \nu_{E}(x))-\chi_{E}(x + \tau \nu_{E}(x)) )\left(1 + H_E(x) \tau  + K_E(x) \tau^2\right)\, d\tau.
\end{equation}
 Then it holds $d_{H^{-1}}(F;E) <\infty$ if and only if
\begin{equation}
\label{eq:zero-condition}
\int_{\Sigma_i} \xi_{F,E} \, d \H^2 = 0 
\end{equation}
for every component $\Sigma_i$ of $\pa E$. 
\end{lemma}

\begin{proof}
The lemma follows from the co-area formula  and the fact that  the tangential Jacobian of $\Phi_\tau$ defined in \eqref{def:mapping-Phi} is 
\[
 J_{\pa E}\Phi_\tau(x) =  1 +  H_E(x) \tau + K_E(x) \tau^2. 
\]
In fat, the assumption $F \Delta E \subset  \mathcal{N}_{r_0}(\pa E)$, combined with the co-area formula and a change of variables, implies that we may write for $f \in H^1(\pa E)$
\[
\begin{split}
\int_{\R^3}   f( \pi_{\pa E}(x))&(\chi_F(x)-\chi_E(x))\,dx = \int_{\mathcal{N}_{r_0}(\pa E)}   f( \pi_{\pa E}(x))(\chi_F(x)-\chi_E(x))\,dx\\
&= \int_{-r_0}^{r_0} \int_{\{x:\, d_E(x) = \tau\}}  f( \pi_{\pa E}(x))(\chi_F(x)-\chi_E(x))\,d \H^2(x)\, d\tau   \\
&=  \int_{-r_0}^{r_0}  \int_{\pa E}  f(x)  (\chi_F(\Phi_\tau(x))-\chi_E(\Phi_\tau(x))) J_{\pa E}\Phi_\tau(x) \, d\H^{2}(x) \, d\tau \\
&=\int_{\pa E}  f(x) \xi_{F,E}(x)  \, d\H^{2}(x) . 
\end{split}
\]
If \eqref{eq:zero-condition} holds, we may assume that   also $f$ satisfies $\int_{\Sigma_i} f \, d \H^2 = 0$ in every component $\Sigma_i$. By Poincar\'e inequality it holds
\[
\int_{\pa E}  f(x) \xi_{F,E}(x)  \, d\H^{2}(x) \leq C \|f\|_{L^2(\pa E)} \leq C \|\nabla_{\pa E} f\|_{L^2(\pa E)}
\]
and therefore  $d_{H^{-1}}(F;E)$ is finite.  It is also immediate that \eqref{eq:zero-condition} is a necessary condition for  $d_{H^{-1}}(F;E) <\infty$. 
\end{proof}

 If in addition to the assumptions of Lemma \ref{lem:distance-finite}, the set $F$ is $C^1$-regular and $\pa F$ is a normal deformation of  $\pa E$, i.e., 
\[
\pa F = \{ x + \psi(x) \nu_E(x) : x \in \pa E \},
\]
 then we may write  the function \eqref{def:xi-original} explicitly as
\begin{equation}\label{def:normal-velo}
  \xi_{F,E} = \psi + \frac{H_E}{2} \psi^2  + \frac{K_E}{3}\psi^3.
\end{equation}
We may also characterize the function $f \in H^1(\pa E)$ attaining the supremum in \eqref{def:akka-meno-uno} as a solution of the equation
\begin{equation}\label{eq:akka-meno-uno}
-\Delta_{\pa E} f =\frac{\xi_{F,E}}{d_{H^{-1}}(F;E)} \qquad \text{on } \pa E. 
\end{equation}
Finally note that if $E$ is in addition  $C^{2,\alpha}$-regular, then $f$ in \eqref{eq:akka-meno-uno} is also $C^{2,\alpha}$-regular. This follows from standard elliptic regularity estimates.

Let us define the standard $H^{-1}$-norm for functions on  $\pa E$ as
\begin{equation}\label{def:akka-meno-uno-norm}
\|f\|_{H^{-1}(\pa E)} = \sup_{\|\nabla_{\pa E} g\|_{L^2(\pa E)}\leq1} \int_{\pa E} f g \, d \H^{2}.
\end{equation}
Then   we have trivially  the following interpolation inequality   
\begin{equation}\label{eq:interpolation-akka-meno-uno}
\|f\|_{L^{2}(\pa E)} \leq  \|\nabla_{\pa E } f\|_{L^{2}(\pa E)}^{\frac12} \|f\|_{H^{-1}(\pa E)}^{\frac12} 
\end{equation}
for all $f \in H^1(\pa E)$. If $f$ has zero average on every components of $\Sigma_i$, then Poincar\'e inequality implies $\|f\|_{H^{-1}(\pa E)}  \leq C \|f\|_{L^2(\pa E)} $.   Moreover, by the proof of Lemma \ref{lem:distance-finite} it holds 
\begin{equation}\label{eq:norm-equivalence}
d_{H^{-1}}(F;E)= \|\xi_{F,E} \|_{H^{-1}(\pa E)},  
\end{equation}
where $\xi_{F,E}$ is defined in \eqref{def:xi-original}.

\subsection{Useful formulas}

In case  $\pa F$ is a normal deformation of $\pa E$ and is parametrized by  the heightfunction $\psi$, we may write the mean curvature $H_F$ in terms of $\psi$.
These calculations are standard and can be found e.g. in \cite{MantegazzaBook} (see also \cite{FJM3D}). However, we need this formula with an explicit dependence on the reference set $\pa E$, which we did not find in the literature. 
The key point is that the expansion of the mean curvature of $F$ depends on the derivative of 
$B_E$, which presents a significant challenge in the forthcoming analysis. This dependence leads to a loss of regularity when expanding the mean curvature using the height function. Specifically, it requires more regularity from $E$ than what the formula provides for $F$. We address this issue by using the fact that the term involving the derivatives of $B_E$ is quadratically small.

We assume that $E$ is $C^3$-regular, satisfies the UBC with radius $r$, and $F$ is $C^2$-regular such that $F \Delta E \subset  \mathcal{N}_r(\pa E)$ and we may write its boundary as 
\begin{equation}
\label{def:graph-of-E}
\pa F = \{ x + \psi(x) \nu_E(x) : x \in \pa E \}
\end{equation}
for $\psi \in C^2(\pa E)$. We define a diffeomorphism $\Psi: \pa E \to \pa F$ as
\begin{equation}
\label{def:diffeo-Psi}
\Psi(x) := x + \psi(x) \nu_E(x) .
\end{equation} 

Let us first make some easy observations. Since  $F \Delta E \subset  \mathcal{N}_r(\pa E)$, the projection $\pi_{\pa E}$ is well defined on $\pa F$ and it holds  for $y \in \pa F$, setting $x = \pi_{\pa E}(y)$,
\begin{equation}
\label{eq:differential-proj}
d_E(y) = \psi(x), \qquad  \nabla \pi_{\pa E}(y) =  I - \nu_E(x) \otimes \nu_E(x)  - \psi(x)  \, B_E(x) \big(I+\psi(x) B_E(x)\big)^{-1}, 
\end{equation}
where the latter follows from \eqref{eq:derivative-proj}. We may extend a given function $f \in C^1(\pa E)$ as $f \circ \pi_{\pa E} : \mathcal{N}_r(\pa E) \to \R$. In particular, this function is defined on $\pa F$ and  we may  calculate the tangential gradient of $f \circ \pi_{\pa E} : \pa F \to \R$   at $y \in \pa F$  as
\begin{equation}
\label{eq:change-variables-function}
\nabla_{\pa F} (f \circ \pi_{\pa E})(y) =  (\nabla_{\pa F} \pi_{\pa E}(y))^T \nabla_{\pa E}f(\pi_{\pa E}(y)).
\end{equation}
Indeed, we obtain \eqref{eq:change-variables-function} by using the fact  that in  $\mathcal{N}_r(\pa E)$  it holds $f \circ \pi_{\pa E} = f \circ \pi_{\pa E} \circ \pi_{\pa E}$.
We differentiate this and have for all $y \in \mathcal{N}_r(\pa E)$
\[
\nabla (f \circ \pi_{\pa E})(y) =  \big(\nabla  \pi_{\pa E}(y)\big)^T \nabla  (f \circ \pi_{\pa E})(\pi_{\pa E}(y)).
\]
The equality \eqref{eq:change-variables-function} then follows from above, by observing that since $\pi_{\pa E}(y) \in \pa E$ it holds by the definition of tangential gradient 
$\nabla  (f \circ \pi_{\pa E})(\pi_{\pa E}(y)) =  \nabla_{\pa E}  f(\pi_{\pa E}(y))$. 

We may now derive our formula for the expansion of the mean curvature. 
\begin{lemma}
\label{lem:expand-mean-curv}
Assume $E\subset\R^n$ is $C^3$-regular and satisfies the UBC with radius $r$, and $F\subset\R^n$ is $C^2$-regular such that $F \Delta E  \subset \mathcal{N}_r(\pa E)$ and we may write its boundary as 
\eqref{def:graph-of-E}  for $\psi \in C^2(\pa E)$. Then it holds for $x \in \pa E= \Sigma$
\[
H_F(\Psi(x)) = - \Delta_{\pa E}\psi(x)  + H_E(x) + R_0(x),
\]
where $\Psi$ is defined in \eqref{def:diffeo-Psi}. The error term is of the form 
\begin{equation}
\label{eq:error-R-0}
R_0 = \la  A_1(\psi B_\Sigma ,\bar \nabla \psi ),  \bar \nabla^2 \psi  \ra +  \la  A_2(\psi B_\Sigma,\bar \nabla \psi ) ,  \bar \nabla (\psi B_\Sigma)\ra  + a_0(\psi,\bar \nabla \psi , B_\Sigma)
\end{equation}
where $A_1,A_2$  are smooth tensor fields  such that $A_i(0,0) = 0$ and $a_0$ is a smooth function with $a_0(0,0,\cdot) = 0$.
\end{lemma}
\begin{proof}
We observe that for $\Psi : \pa E \to \pa F$  defined as in \eqref{def:diffeo-Psi} it holds  $ \Psi^{-1}= \pi_{\pa E} : \pa F \to \pa E$. 
For every tangent vector  $ \tau \in T_x \pa E$  it holds $\nabla_{\pa E} \Psi(x) \tau \in  T_{\Psi(x)} \pa F$  and 
\begin{equation}\label{eq:meancurva-exp-00}
\nabla_{\pa E} \Psi(x) \tau = (I+ \psi(x)B_E(x))\tau+ (\nabla_{\pa E} \psi \cdot \tau)  \nu_E(x).
\end{equation}
We define $N : \pa E \to \R^n$ as 
\begin{equation}\label{eq:meancurva-exp-0}
N(x)= -(I+\psi(x) B_E(x))^{-1} \nabla_{\pa E} \psi(x) + \nu_{ E}(x). 
\end{equation}
Then it holds   for every  $ \tau \in T_x \pa E$ that 
\[
\big(\nabla_{\pa E} \Psi(x) \tau \big) \cdot N(x) =  0. 
\]
It follows from \eqref{def:graph-of-E} that $\nu_F(\Psi(x)) \cdot \nu_E(x) \geq 0$. Therefore  at every  $y \in \pa F$ the vector $N(\pi_{\pa E}(y))$ is in the direction of  the normal $\nu_F(y)$, i.e., 
\begin{equation}\label{eq:meancurva-exp-1}
	\nu_F(y)= \frac{N (\pi_{\pa E}(y))}{\vert N (\pi_{\pa E}(y)) \vert}.
\end{equation}
Therefore it holds 
\begin{equation}\label{eq:meancurva-exp-2}
	\begin{split}
H_F(y) &= (\div_{\pa F}\nu_F)(y) \\
&= \frac{1}{\vert N (\pi_{\pa E}(y)) \vert}  (\div_{\pa F} (N \circ \pi_{\pa E}))(y)  + \nabla_{\pa F} \left( \vert N (\pi_{\pa E}(y)) \vert^{-1} \right) \cdot N (\pi_{\pa E}(y))\\
&=  \frac{1}{\vert N (\pi_{\pa E}(y)) \vert}  (\div_{\pa F} (N \circ \pi_{\pa E}))(y),
\end{split}
\end{equation}
where the  last equality follows from the fact that $N(\pi_{\pa E}(y))$ is in the direction of $\nu_F(y)$.

We split $N = N_1 + N_2$ as 
\[
N_1 = -(I+\psi B_E)^{-1} \nabla_{\pa E} \psi  \quad \text{and} \quad N_2  = \nu_{ E}. 
\]
Note that $N_1$ and $N_2$ are orthogonal and $N_1$ is a tangent vector field, by which we mean that  $N_1(x) \in T_x \pa E$ for every $x \in \pa E$.  By \eqref{eq:meancurva-exp-2} we need to calculate $\div_{\pa F} (N_1 \circ \pi_{\pa E})$ and $\div_{\pa F} (N_2 \circ \pi_{\pa E})$. 
We  write 
\begin{equation}\label{eq:meancurva-exp-3}
\vert N (\pi_{\pa E}(y)) \vert = 1 + a_1(\psi B_\Sigma,\bar \nabla \psi ),
\end{equation}
where the function satisfies  $a_1(0,0) = 0$. 

Let us first calculate $\div_{\pa F} (N_1 \circ \pi_{\pa E})$. Fix $y \in \pa F$ and denote $x = \pi_{\pa E}(y)$.  We use \eqref{eq:change-variables-function} component-wise to infer
\begin{equation}\label{eq:meancurva-exp-4}
\div_{\pa F} (N_1 \circ \pi_{\pa E})(y)  =  \mathrm{Tr} (\nabla_{\pa E} N_1(x) \nabla \pi_{\pa E}(y) P_{\pa F}(y)),
\end{equation}
where $P_{\pa F} = I - \nu_F \otimes \nu_F$. By \eqref{eq:meancurva-exp-1}, \eqref{eq:meancurva-exp-3} and by the definition of $N$ in \eqref{eq:meancurva-exp-0}  we find that there are smooth tangent vector fields $\hat A_i = \hat A_i (\psi B_E,\nabla_{\pa E} \psi)$, with $i =1,2$, and a smooth  function $a_2 = a_2 (\psi(x) B_E(x),\nabla_{\pa E} \psi(x))$ such that   
\begin{equation}\label{eq:meancurva-exp-5}
P_{\pa F}(y) - P_{\pa E}(x) = -\hat A_1 \otimes \hat A_1 +  \hat A_2 \otimes \nu_E(x) +  \nu_E(x) \otimes \hat A_2 + a_2 \, \nu_E(x) \otimes \nu_E(x)  ,
\end{equation}
and it holds  $\hat A_i(0,0) = 0$, $a_2(0,0) = 0$. We may thus write   \eqref{eq:meancurva-exp-4} using \eqref{eq:differential-proj},  \eqref{eq:meancurva-exp-5} and recalling that $N_1 = -(I+\psi B_E)^{-1} \nabla_{\pa E}\psi $ is a tangent vector field
\begin{equation}\label{eq:meancurva-exp-6}
\div_{\pa F} (N_1 \circ \pi_{\pa E}) =  - \Delta_{\pa E}\psi +  \la  \tilde A_1(\psi B_\Sigma ,\bar \nabla \psi ),  \bar \nabla^2 \psi  \ra +  \la \tilde  A_2(\psi B_\Sigma,\bar \nabla \psi ) ,  \bar \nabla (\psi B_\Sigma)\ra ,
\end{equation}
where $\tilde A_1,\tilde A_2$  are smooth tensor fields  such that $\tilde A_i(0,0) = 0$. 

Let us then calculate  $\div_{\pa F} (N_2 \circ \pi_{\pa E})$.  We have $N_2 \circ \pi_{\pa E} = \nu_E \circ \pi_{\pa E}  = \nabla d_E$. By \eqref{eq:Hessian-distance} it holds
\[
\nabla^2d_E(y)=B_E(x)(I+\psi(x)B_E(x) )^{-1}.  
\]
Therefore we have by \eqref{eq:meancurva-exp-5} (denote the vector field  by $\hat A_1= \hat A_1(\psi B_E,\nabla_{\pa E} \psi)$) for $y \in \pa F$ that 
\begin{equation}\label{eq:meancurva-exp-7}
	\begin{split}
		\div_{\pa F}(N_2 \circ \pi_{\pa E}) &= \mathrm{Tr}(\nabla_{\pa F} \nabla d_E)  =  \mathrm{Tr}\big(P_{\pa F}(y) \nabla^2 d_E(y)) \big)  \\
		&= \mathrm{Tr}(\nabla^2 d_E(y)) -(\nabla^2 d_E(y) \, \hat A_1) \cdot \hat A_1\\
		&=  \mathrm{Tr}\big(B_E(x)(I+\psi(x)B_E(x) )^{-1})\big) - \big(B_E(x)(I+\psi(x)B_E(x) )^{-1}  \, \hat A_1\big) \cdot \hat A_1\\
		&= H_E(x) +  a_3(\psi, \nabla_{\pa E} \psi, B_E), 
	\end{split}
\end{equation}
where $a_3$ is a smooth function such that  $a_3(0,0,\cdot )= 0$. The claim then follows from \eqref{eq:meancurva-exp-2}, \eqref{eq:meancurva-exp-3},  from $N = N_1 + N_2$, \eqref{eq:meancurva-exp-6}  and from \eqref{eq:meancurva-exp-7}. 
\end{proof}

We remark that if $E$ and $F$ are as in Lemma \ref{lem:expand-mean-curv}, then we may use \eqref{eq:meancurva-exp-5} to write the tangential derivative of a function $f : \pa E \to \R$ in \eqref{eq:change-variables-function} on $y \in \pa F$, setting $x = \pi_{\pa E}(y)$, as 
\begin{equation}
\label{eq:change-variables-function-2}
\nabla_{\pa F} (f \circ \pi_{\pa E})(y) = \big(I + A(\psi(x), \nabla_{\pa E} \psi(x), \nu_E(x),B_E(x))\big)\nabla_{\pa E}f(x)
\end{equation}
for a smooth matrix-valued function with $A(0,0,\cdot, \cdot) = 0$.

\subsection{Almost minimizers of the perimeter}

In this section we discuss  properties of almost minimizers of the perimeter.  We will  use the following definition. 
\begin{definition}
\label{def:lambda-mini}
A set of finite perimeter $F\subset \R^n$ is a $(\Lambda,\alpha)$-minimizer of the perimeter for $\Lambda \geq 0$ and $\alpha \in [0,\tfrac{1}{n})$ if for  every set of finite perimeter   $G \subset \R^n$  it holds
\[
P(F) \leq P(G) + \Lambda |F\Delta G|^{1-\alpha}. 
\]
Moreover, we say that  $F$ is  a $\Lambda$-minimizer of the perimeter if it is  $(\Lambda,0)$-minimizer of the perimeter.
\end{definition}
We point out that the definition of  $(\Lambda,\alpha)$-minimizer is closely related to  the one of  $\omega$-minimizer of the perimeter (see \cite{Tamanini1984}), where  a set of finite perimeter $F\subset \R^n$ is an  $(\omega_0,r_0,\beta)$-minimizer of the perimeter if for all $x$ and $r \leq r_0$ and  every set of finite perimeter   $G \subset \R^n$ with $G\Delta F \subset B_r(x)$ it holds
\[
P(F) \leq P(G) + \omega_0 r^{n-1+\beta}. 
\]
It is clear that if a set is a $(\Lambda,\alpha)$-minimizer of the perimeter  according to  Definition \ref{def:lambda-mini}, then it is  $(\omega_0,r_0,\beta)$-minimizer for $\omega_0 = \Lambda |B_1|^{1-\alpha}$, say for $r_0 = 1$, and $\beta = 1 -n \alpha$. 
Therefore  Definition \ref{def:lambda-mini} is stronger than the $\omega$-minimality defined above as it does not require the  perturbation to be local. 

It is well-known that  if  $E  \subset \R^n$ satisfies the UBC with radius $r_0$, then it is a $\Lambda$-minimizer of the perimeter, where $\Lambda = \Lambda(n,r_0)$. We state this result below as a lemma whose proof can be found   in  \cite[Lemma 4.1]{AFM}.

\begin{lemma}
\label{lem:UBC-lambda-min}
Assume that $E \subset \R^n$ is bounded and satisfies the UBC with  radius $r_0$, then it is  a $\Lambda$-minimizer of the perimeter, i.e., every set of finite perimeter   $G \subset \R^n$  it holds
\[
P(E) \leq P(G) +  \Lambda |E\Delta G|, 
\]
for $\Lambda = \frac{C}{r_0}$, for a dimensional constant  $C$. In particular, $P(E) \leq \Lambda|E|$.
\end{lemma}
We remark that if $E \subset \R^n$  satisfies the UBC with radius $r_0$ and $E \subset B_{R_0}$, then we may bound its perimeter  by Lemma \ref{lem:UBC-lambda-min}  as $P(E) \leq \Lambda|E| \leq \Lambda |B_{R_0}|$. 

The following lemma states some useful regularity properties that almost minimizers of the perimeter introduced in Definition \ref{def:lambda-mini} inherit from perimeter minimizers \cite{Tamanini1984}.

\begin{lemma}
\label{lem:omega-mini}
Assume that $E \subset B_{R_0}\subset\R^n$ satisfies the UBC with radius $r_0$  and is uniformly $C^{2,\beta}$-regular with constants $r_0, C_0$, and $F$ is  $(\Lambda, \alpha)$-minimizer of the perimeter   for $\alpha \in (0,\tfrac{1}{n})$. Then for  $\gamma < \min\{\beta, \frac12 - \frac{n \alpha}{2} \}$ there is  $\delta_0 \in (0,r_0/2)$, depending on $R_0,r_0,C_0,\Lambda, \alpha, \gamma$ and the dimension, such that 
if 
\[
F \Delta E \subset \mathcal{N}_{\delta_0}(\pa E)
\]
then there is $\psi \in C^{1,\gamma}(\pa E)$ such that 
\[
\pa F = \{ x + \psi(x)\nu_E(x) : x \in \pa E\}.
\]
Moreover, for every $\eps>0$ there is $\delta_0 = \delta_0(\eps)$ such that  $\|\psi\|_{C^{1,\gamma}(\pa E)} \leq \eps$.
\end{lemma} 

\begin{proof}
We first note that by the assumption $F \Delta E \subset \mathcal{N}_{\delta_0}(\pa E)$  for every $x \in \pa E$ the set  $ \pa F \cap B_{\delta_0}(x)$ is non-empty. 

 Let us fix $\eps>0$. We claim that there is $\delta_0= \delta_0(\eps) \in (0, \frac{r_0}{80}\eps)$ such that    if  $F$ and $E$ are as in the assumptions,  then  for every  $x \in \pa E$   and every $y \in \pa F \cap B_{\delta_0}(x)$ it holds
\begin{equation}
\label{eq:normal-continuous}
|\nu_E(x) - \nu_F(y)| < \eps. 
\end{equation}
  We argue  by contradiction and assume that there exist $E_k, F_k$  satisfying the assumptions with $E$ and $F$ replaced by $E_k$ and $F_k$, respectively, such that $F_k \Delta E_k \subset \mathcal{N}_{\frac{1}{k}}(\pa E_k)$ and points $x_k \in \pa E_k$, $y_k \in \pa F_k \cap B_{\frac{1}{k}}(x_k)$ such that 
\begin{equation}\label{eq:contr-normal-continuous}
|\nu_{E_k}(x_k) - \nu_{F_k}(y_k)| \geq \eps
\end{equation}
for some $\eps>0$. By passing to a subsequence we may assume that  $x_k \to x$, $y_k \to x$,  $E_k \to E$ in $C^{2,\beta'}$ for all $0<\beta'<\beta$ and $F_k\to E$ in the Hausdorff distance, $E$  satisfies the UBC and  is  in particular  a $\Lambda$-minimizer of the perimeter. Then $ \nu_{E_k}(x_k) \to \nu_{E}(x)$, while  from  $(\Lambda, \alpha)$-minimality of the sets $F_k$ it follows by \cite[Section 1.9]{Tamanini1984}  that  $\nu_{F_k}(y_k) \to \nu_{E}(x)$, which contradicts \eqref{eq:contr-normal-continuous}.

The lemma follows from \eqref{eq:normal-continuous} by standard regularity argument. Indeed, let us fix $x_0 \in \pa E$. We may assume $x_0 = 0$ and $\nu_E(0) = e_n$. Since $E$ satisfies the UBC with radius $r_0$,    we may write $ E \cap B_{r_0/2}$ as a subgraph of a function $f: B_{r_0/2}^{n-1} \to \R$ with $\|f\|_{C^{1,1}(B_{r_0/2}^{n-1})} \leq \frac{10}{r_0}$, assuming $r_0 \leq 1$, see \eqref{eq:UBC-graph-function}.  Then it holds  $|\nu_E(x) - e_n | < \eps$   for all $x \in \pa E \cap B_{r_\eps}$, where $r_\eps = \frac{r_0}{20}\eps$. We note that  $\delta_0 < \frac{r_0}{80}\eps$ implies  $\delta_0 < \frac{r_\eps}{4}$.   Then by  \eqref{eq:normal-continuous} we have $|\nu_F(y) - e_n | < 2\eps$ for all $y \in \pa F \cap B_{\frac{3r_\eps}{4}}$. Choose any point  $y_0 \in \pa F \cap B_{\delta_0}$. By the previous inequality and the perimeter density estimates  for $(\Lambda,\alpha)$-minimizers we then deduce that the excess 
\[
\mathcal{E}\Big(y_0,\frac{r_\eps}{2}\Big) =\min_{|\omega| = 1 }  \frac{1}{(r_\eps/2)^{n-1}}  \int_{\pa F \cap B_{\frac{r_\eps}{2}}(y_0)} |\nu_F(y) - \omega |^2 \, d \H^{n-1}(y)\leq C\eps^2,
\]
provided $r_\eps<r_1=r_1(n,\Lambda,\alpha)$, for some constant $C=C(n,\Lambda,\alpha)$. Then, 
 it follows from the so-called $\eps$-regularity theorem, see for instance \cite{Tamanini1984}, and from the inclusion $B_{r_\eps/4} \subset B_{r_\eps/2}(y_0) $   that there is $\vphi : B_{r_\eps/4}^{n-1} \to \R$ such that 
\[
 F \cap B_{r_\eps/4} = \{ (y',y_n) \in \R^n : y_n < \varphi(y') \} \cap  B_{r_\eps/4} 
\]
 with   $\|\vphi\|_{C^{1,\gamma'}(B_{r_\eps/4}^{n-1})}\leq C$ and $\gamma' =  \frac12 -\frac{ n\alpha}{2}$. The existence of the heightfunction $\psi \in C^{1,\gamma}(\pa E)$  then follows from  the fact that $\pa E$ is uniformly $C^{2,\beta}$-regular, see \cite[Section 1.2]{Antonia} for the precise argument. 
 The smallness of the norm $\|\psi\|_{C^{1,\gamma}(\pa E)}$ for $\gamma < \min\{ \beta, \gamma'\}$ when $\delta_0$ is small, follows from interpolation inequality \eqref{eq:Holder-interpolate} observing that $\|\psi\|_{C^0} \leq \delta_0$. 
\end{proof}

As stated in Lemma \ref{lem:UBC-lambda-min},  if  $E \subset B_{R_0}$ satisfies the UBC with radius $r_0$, then it is a $\Lambda$-minimizer of the perimeter. Moreover for every $s \in (-\frac{r_0}{2}, \frac{r_0}{2})$ also the set 
\begin{equation}
\label{eq:minkowski_sum}
E_s = \{ x \in \R^n : d_E(x) <s \}
\end{equation}
satisfies the UBC with  radius $r_0/2$. Thus Lemma \ref{lem:UBC-lambda-min} implies that if $G \subset \R^n$ is a set of finite perimeter, then for every $s \in (-\frac{r_0}{2}, \frac{r_0}{2})$  it holds 
\begin{equation}
\label{eq:minkowski_ineq}
P(G \cap E_s) \leq  P(G) + \Lambda |G \setminus E_s| \quad \text{and} \quad  P(G \cup E_s) \leq P(G) + \Lambda |E_s \setminus G|,
\end{equation}
for $\Lambda = \Lambda (r_0,n)$, by possibly increasing $\Lambda$.  

We will consider constrained minimization problems and to this aim we introduce the following definition.
\begin{definition}\label{def:constr-lambda-mini}
A set of finite perimeter $F\subset\R^n$ is a constrained $(\Lambda,\alpha)$-minimizer of the perimeter for $\Lambda>0$ and $\alpha\in[0,\frac1n)$ if for every set of finite perimeter $G$ such that 
\[
G \Delta F \subset \overline{\mathcal{N}_{\delta}(\pa E)} \quad  \text{and } \quad |G|= |F| 
\]
it holds
\[
P(F) \leq P(G) + \Lambda |F\Delta G|^{1-\alpha} .
\]
\end{definition}

\begin{lemma}
\label{lem:constraint-omega-mini}
Assume that $E \subset \R^n$ satisfies the UBC with radius $r_0$, $E \subset B_{R_0}$ and let $F$ be a constrained $(\Lambda,\alpha)$-minimizer such that  $|F| = |E|$ and 
\[
F \Delta E \subset \overline{\mathcal{N}_{\delta}(\pa E)}
\]
for $\delta \in (0,r_0/4)$. Then $F$ is a $(\Lambda_1,\alpha)$-minimizer of the perimeter (without constraint) according to Definition \ref{def:lambda-mini}, where 
$\Lambda_1=\Lambda_1(\Lambda, r_0, R_0,n)$.   
\end{lemma}

\begin{proof}
Let us fix a set of finite perimeter $G \subset \R^n$.  We first notice that we may assume $ |F\Delta G| \leq 1$. Indeed, if not then the inequality $P(F) \leq P(G) + \Lambda_1 |F\Delta G|^{1-\alpha}$ holds trivially by choosing $\Lambda_1 \geq P(F)$. Note also that  we may bound the perimeter of $F$ by choosing $G= E$ and using the assumption on $F$ to have 
\[
P(F) \leq P(E) +  \Lambda |F\Delta E|^{1-\alpha} \leq P(E) + \Lambda|B_{R_0}|.
\]
The above estimate leads to the bound on the perimeter of $F$ recalling that Lemma \ref{lem:UBC-lambda-min} yields $P(E) \leq \Lambda |B_{R_0}|$.

The rest of the proof is based on a calibration  argument and to this aim we recall our notation \eqref{eq:minkowski_sum} for $E_s = \{ x \in \R^n : d_E(x) <s \}$ and that $E_s$ satisfies the UBC for radius $r_0/2$  for all  $s \in (-\frac{r_0}{2}, \frac{r_0}{2})$. 
The following argument will only rely on the estimate \eqref{eq:minkowski_ineq} and therefore all  the constants  below,   denoted simply  by $\Lambda_0$,  will depend only on $r_0$ and the dimension $n$.  

Let $G \subset \R^n$ be such that  $ |F\Delta G| \leq 1$. We define 
\[
\tilde G = (G \cap E_{\delta}) \cup  E_{-\delta}. 
\]
Then it holds $ \tilde G \Delta E  \subset \overline{\mathcal{N}_{\delta}(\pa E)}$ and  $\tilde G \Delta  G = (G \setminus E_{\delta}) \cup (E_{-\delta} \setminus G)$. Since $F \Delta E  \subset  \overline{\mathcal{N}_{\delta}(\pa E)}$ then $F \subset \overline E_{\delta}$ and $\text{int} \,  E_{-\delta} \subset  F$, which implies  
\begin{equation}
\label{eq:minkowski_lemma-1}
|G \Delta \tilde  G| \leq  |G \Delta F|.
\end{equation}
Using \eqref{eq:minkowski_ineq}  and \eqref{eq:minkowski_lemma-1} we may estimate the perimeter as
\begin{equation}
\label{eq:minkowski_lemma-2}
P(\tilde G) \leq P(G) + \Lambda_0 |G \Delta  \tilde G| \leq P(G) + \Lambda_0 |G \Delta F|. 
\end{equation}

It holds either $|\tilde G | \geq |F|$ or  $|\tilde G | \leq |F|$. Let us assume the former, the other case being similar. Define a continuous  function $\omega: [-\delta, \delta] \to \R$,
\[
\omega(s) = |E_s \cap \tilde G|. 
\]
It holds by the construction of $\tilde G$ and from $E_{-\delta} \subset F$ that 
\[
\omega(\delta) = |\tilde G| \geq |F| \quad \text{ and  }\quad \omega(-\delta) = |\tilde G \cap E_{-\delta}| =  | E_{-\delta}| \leq |F|.
\]
Therefore there exists $\tau \in [-\delta, \delta] $ such that $\omega(\tau)  = |F|$, and thus the set $G_\tau = E_{\tau} \cap \tilde G$ satisfies $|G_\tau | = |F|$. Note also that it follows from $G_\tau \subset \tilde G$ and \eqref{eq:minkowski_lemma-1} that 
\begin{equation}
\label{eq:minkowski_lemma-3}
\begin{split}
|\tilde G \Delta G_\tau| &= |\tilde G | - |G_\tau| =   |\tilde G | - |F| \leq |\tilde G \Delta F| \\
&\leq |G \Delta F| +|G \Delta \tilde G| \\
&\leq 2  |G \Delta F|.  
\end{split}
\end{equation}
Moreover, using \eqref{eq:minkowski_ineq} and \eqref{eq:minkowski_lemma-3} we have 
\begin{equation}
\label{eq:minkowski_lemma-4}
P(G_\tau) \leq P(\tilde G) + \Lambda_0 |\tilde G \Delta G_\tau| \leq P(\tilde G) + 2\Lambda_0 |G \Delta F|.  
\end{equation}

Since $G_\tau \Delta E \subset \overline{\mathcal{N}_{\delta}(\pa E)}$ and $|G_\tau| = |F|$, we have by the assumption on $F$ that 
\[
P(F) \leq P(G_\tau) + \Lambda |G_\tau \Delta F|^{1-\alpha} \leq P(G_\tau) + \Lambda (|G \Delta F|^{1-\alpha} + |G \Delta \tilde G|^{1-\alpha} + |\tilde G \Delta G_\tau|^{1-\alpha} ) . 
\] 
The claim then follows from  \eqref{eq:minkowski_lemma-1},   \eqref{eq:minkowski_lemma-2},  \eqref{eq:minkowski_lemma-3} and \eqref{eq:minkowski_lemma-4}. 
\end{proof}

\section{Geometric estimate}

 We need a higher order version of the fact  that if a bounded set $E$ satisfies the UBC, then it is a $\Lambda$-minimizer of the perimeter (Lemma \ref {lem:UBC-lambda-min}). While the result in  Lemma \ref {lem:UBC-lambda-min} follows from  standard calibration argument, its higher  order version, which is  formulated  below, is considerably more challenging due to the degenerate nature of the $H^{-1}$-distance  \eqref{def:akka-meno-uno}.   Note that, since the UBC condition implies that $E$ is of class $C^{1,1}$, this allows us to define the Sobolev space $H^2(\Sigma)$, $\Sigma=\pa E$ as discussed in the previous section.  
\begin{proposition}
\label{prop:geometric}
Assume $E \subset \R^3$ satisfies the UBC with radius $r_0$, $E \subset B_{R_0}$, and $\|H_E\|_{H^2(\Sigma)} \leq K_0$. Then there exist $\delta_1 \in (0,r_0)$ and $\Lambda_2 \geq 1$, depending on $r_0, R_0$ and $K_0$, such that  for every set of finite perimeter $F \Delta E  \subset  \mathcal{N}_{\delta_1}(\pa E)$ it holds
\begin{equation}
\label{eq:geometric-1}
P(E) \leq P(F) +\Lambda_2   d_{H^{-1}}(F;E),
\end{equation}
where $ d_{H^{-1}}(F;E)$ is defined in \eqref{def:akka-meno-uno}. Moreover, if $F$ is $C^{1}$-regular and there is  $ \psi \in C^{1}(\partial E)$ such that 
\begin{equation}\label{psi-delta1}
 \partial F = \{ x+ \psi(x)\nu_{E}(x): x \in \partial E\} \quad \text{ with } \quad  \| \psi \|_{C^{1}(\partial E)} <\delta_1 
\end{equation}
then
\begin{equation}
\label{eq:geometric-2}
\frac{1}{4}\|\nabla_{\pa E} \psi \|_{L^2(\pa E)}^2+ P(E) \leq P(F)+ \Lambda_2 d_{H^{-1}}(F;E) .
\end{equation}
\end{proposition}

\begin{proof}
We begin by proving the claim \eqref{eq:geometric-2}. Note first that the assumption $\|H_E\|_{H^2(\Sigma)} \leq K_0$ implies $\|B_\Sigma\|_{H^2(\Sigma)} \leq C$. Indeed, given $x\in\pa E$, up to a rotation we can write $ E \cap B_{r_0/2}(x)$ as in \eqref{eq:UBC-graph-function}, where the function $f\in C^{1,1}(B^2_{r_0/2}(x'))$
satisfies in $B^2_{r_0/2}(x')$ the mean-curvature equation with right-hand-side given by the function $\tilde H(x')=H_E((x',f(x')))\in H^2(B^2_{r_0/2}(x'))$. As a result, by standard elliptic estimates we get that $f$ is bounded in $H^4(B^2_{r_0/4}(x'))$, which in turn, implies that $B_{\Sigma}$ is bounded in $H^2(\Sigma)$ and, by the Sobolev embedding theorem, that $E$ is uniformly $C^{2,\alpha}$-regular for any $\alpha \in (0,1)$ and it holds 
\begin{equation}
\label{eq:geometric-curv-bound}
\|B_E\|_{L^\infty} + \|\bar \nabla B_\Sigma\|_{L^p} \leq C_p,
\end{equation} 
where $C_p = C_p(r_0,R_0,  K_0,p)$ for every $p \in (1,\infty)$.

Since $E$ is uniformly  $C^{2,\alpha}$-regular we  may  define a diffeomorphism $\Psi : \pa E \to \pa F$,
\[
\Psi(x) =  x+ \psi(x)\nu_{E}(x). 
\]
Let us fix a point $x \in \pa E$ and choose   an orthonormal  basis $\tau_1,\tau_2$ of the tangent hyperplane $T_x\pa E$ such that $B_E(x) \tau_i = \kappa_i(x) \tau_i$,  where $\kappa_i(x)$ are the principal curvatures of $\pa E$. By \eqref{eq:meancurva-exp-00} it holds
\[
\nabla_{\pa E} \Psi(x) \tau_i = (1+ \psi(x) \kappa_i(x)) \tau_i + \bar \nabla_i \psi(x)\, \nu_E(x), 
\]
where $\bar \nabla_i \psi(x)$ denotes the derivative in direction $\tau_i$.  Hence,  the tangential Jacobian is 
\begin{equation}
\label{eq:jacobian-psi}
J_{\pa E} \Psi = \left( (1+ \psi \kappa_1)^2(1+ \psi \kappa_2)^2 + (1+ \psi \kappa_1)^2 |\bar \nabla_2 \psi |^2 + (1+ \psi \kappa_2)^2 |\bar \nabla_1 \psi |^2  \right)^{\frac12}.
\end{equation}
 Since $\sqrt{1 +\tau} \geq 1 + \frac{\tau}{2} - \tau^2$ for $|\tau| \leq \frac12$, we deduce that 
\[
J_{\pa E} \Psi \geq 1  + H_E \psi + \frac{1}{3} |\bar \nabla \psi|^2 - C\psi^2
\]
when $\|\psi \|_{C^1} \leq \delta_1$ for $\delta_1$ is small enough. We use  the area formula to conclude
\begin{equation}
\label{eq:geometric-4}
P(F) = \int_{\pa E} J_{\pa E} \Psi \, d \H^2   \geq P(E) + \int_{\pa E} H_E\, \psi \, d \H^2 + \frac13 \| \nabla_{\pa E} \psi\|_{L^2(\pa E)}^2 - C  \| \psi\|_{L^2(\pa E)}^2.
\end{equation}

Next we recall \eqref{eq:norm-equivalence}, i.e.,  $d_{H^{-1}}(F;E)  =  \|\xi_{F,E} \|_{H^{-1}(\pa E) } $,  where $\xi_{F,E}$ is defined in  \eqref{def:normal-velo}.  In the following  we write simply $\xi = \xi_{F,E}$.  We claim that when $\|\psi\|_{C^1}\leq \delta_1$, and when $\delta_1 = \delta(r_0,R_0, K_0)$ is small  enough it holds 
\begin{equation}
\label{eq:geometric-44}
\begin{split}
\frac12|\psi(x)| &\leq |\xi(x)|\leq  2|\psi(x)|,\\
 \|\nabla_{\pa E} \xi\|_{L^2(\pa E)} &\leq 2\| \nabla_{\pa E} \psi \|_{L^2(\pa E)} +  \sqrt{\delta_1}  \|\psi \|_{L^2(\pa E)}. 
\end{split}
\end{equation}
 The  inequalities on the first line follow immediately from \eqref{eq:geometric-curv-bound}. For the inequality on the second line in \eqref{eq:geometric-44} we have by straightforward calculation and by applying \eqref{eq:geometric-curv-bound}
\[
|\nabla_{\pa E} \xi(x)| \leq 2 |\nabla_{\pa E} \psi(x) | +C \psi^2(x) |\bar \nabla B_\Sigma(x)| \qquad \text{for a.e. } x \in \pa E. 
\]
By H\"older's inequality and \eqref{eq:geometric-curv-bound}  we have 
\[
\begin{split}
\|\nabla_{\pa E} \xi\|_{L^2(\pa E)} &\leq 2\| \nabla_{\pa E} \psi \|_{L^2(\pa E)}  + C\|\psi\|_{L^8(\pa E)}^2 \| \bar \nabla B_\Sigma\|_{L^4(\Sigma)}\\ 
&\leq 2\| \nabla_{\pa E} \psi \|_{L^2(\pa E)}  + C\|\psi\|_{L^8(\pa E)}^2.
\end{split}
\]
The claim \eqref{eq:geometric-44} then follows from interpolation inequality in  Proposition \ref{prop:interpolation}, as
\[
C\|\psi\|_{L^8(\pa E)}^2 \leq C\|\psi\|_{W^{1,4}(\pa E)}\|\psi\|_{L^2(\pa E)} \leq C \delta_1 \|\psi\|_{L^2(\pa E)} \leq \sqrt{ \delta_1} \|\psi\|_{L^2(\pa E)}, 
\] 
when $\delta_1$ is small enough. 

 We proceed by using  \eqref{eq:geometric-44}  to write \eqref{eq:geometric-4}    as
\begin{equation}
\label{eq:geometric-45}
P(F)   \geq P(E) + \int_{\pa E} H_E\, \xi \, d \H^2 + \frac13 \| \nabla_{\pa E} \psi\|_{L^2(\pa E)}^2 - C  \| \xi \|_{L^2(\pa E)}^2.
\end{equation}
By the assumptions on $E$ and by \eqref{eq:norm-equivalence}   it holds 
\begin{equation}
\label{eq:geometric-46}
 -\int_{\pa E} H_E\, \xi \, d \H^2 \leq  C \|\xi \|_{H^{-1}(\pa E)} = C d_{H^{-1}}(F;E), 
\end{equation}
where the  constant $C$  depends on $r_0$ and $K_0$. 

By the interpolation inequality \eqref{eq:interpolation-akka-meno-uno} and by \eqref{eq:geometric-44} we have 
\[
\begin{split}
 \| \xi\|_{L^2(\pa E)}^2 &\leq  \|\nabla_{\pa E} \xi\|_{L^2(\pa E)}\|\xi \|_{H^{-1}(\pa E)}\\
&\leq  2\|\nabla_{\pa E} \psi\|_{L^2(\pa E)}\|\xi \|_{H^{-1}(\pa E)} +   \sqrt{\delta_1}  \|\psi \|_{L^2(\pa E)} \|\xi \|_{H^{-1}(\pa E)}\\
&\leq   2\|\nabla_{\pa E} \psi\|_{L^2(\pa E)}\|\xi \|_{H^{-1}(\pa E)} + C\sqrt{\delta_1}  \|\xi \|_{L^2(\pa E)}^2,
\end{split}
\] 
where in the last inequality we used  $ \|\xi\|_{H^{-1}(\pa E)} \leq  C\|\xi \|_{L^{2}(\pa E)}$ and $|\psi| \leq 2 |\xi|$. Therefore when $\delta_1$ is small enough we have 
\[
 \| \xi\|_{L^2(\pa E)}^2 \leq  4\|\nabla_{\pa E} \psi\|_{L^2(\pa E)}\|\xi \|_{H^{-1}(\pa E)} .
\]
It also  holds 
\[
\|\xi \|_{H^{-1}(\pa E)} \leq C\|\xi \|_{L^{2}(\pa E)} \leq C\|\psi \|_{L^{2}(\pa E)}\leq  C \delta_1. 
\]
Therefore when $\delta_1$ is small, we have by the two above inequalities and by \eqref{eq:norm-equivalence} that 
\begin{equation}
\label{eq:geometric-47}
\begin{split}
 \| \xi\|_{L^2(\pa E)}^2 &\leq  4\|\nabla_{\pa E} \psi\|_{L^2(\pa E)}\|\xi \|_{H^{-1}(\pa E)} \\
&\leq \eps \|\nabla_{\pa E} \psi\|_{L^2(\pa E)}^2 + C_\eps \|\xi \|_{H^{-1}(\pa E)}^2\\
&\leq \eps \|\nabla_{\pa E} \psi\|_{L^2(\pa E)}^2 + C_\eps \delta_1 d_{H^{-1}}(F;E) .
\end{split}
\end{equation}
 Therefore we deduce from   \eqref{eq:geometric-45}, \eqref{eq:geometric-46} and  \eqref{eq:geometric-47}  that 
\[
\left(\frac{1}{3} -\eps\right) \|\nabla_{\pa E} \psi \|_{L^2(\pa E)}^2+ P(E) \leq P(F)+C d_{H^{-1}}(F;E).
\]
 Thus we have  \eqref{eq:geometric-2} when $\eps$ is chosen small enough.

\medskip

Let us then prove the first claim \eqref{eq:geometric-1}. We define the functional 
\begin{equation}
\label{def:pena-functional}
 \mathcal{J}(F):= P(F)+ \Lambda_2 d_{H^{-1}}(F;E), \qquad F \Delta E \subset  \overline{\mathcal{N}_{\delta}}(\pa E), 
\end{equation}
where  $\Lambda_2$ is the constant  from \eqref{eq:geometric-2} and  $\delta \in (0,\tfrac{r_0}{4})$. Using direct method of calculus of variations we deduce that $\mathcal{J}$ has a minimizer, say $F$. The first claim  \eqref{eq:geometric-1} follows once we show that $E$ is a minimizer of \eqref{def:pena-functional}, when  $\delta$ is chosen small enough. 

We first claim that the minimizer  $F$ is  a $(\Lambda, \alpha)$-minimizer of the perimeter, for any $\alpha \in (0,\tfrac{1}{3})$, according to Definition \ref{def:lambda-mini}, where $\Lambda$ depends on $r_0, R_0,K_0, \Lambda_2$ and $\alpha$. 
Since for $F$ it clearly holds $ d_{H^{-1}}(F;E)  <\infty$, then by Lemma \ref{lem:distance-finite} we may divide the set  $F$ in pieces $F_i$  according to the components $\Sigma_i$  of $\pa E$, i.e.,  $F_i \Delta E \subset   \overline{\mathcal{N}_{\delta}}(\Sigma_i )$. By taking variations of  each piece $F_i$ separately, we may 
reduce to  compare $F$ with sets  $G \subset \overline{\mathcal{N}_{\delta}}(E) $ such that $G \Delta F \subset \overline{\mathcal{N}_{\delta}}(\Sigma_i)$. For these sets we have by Lemma   \ref{lem:distance-finite}  that $ d_{H^{-1}}(G;E)  <\infty$ if  $|G| = |E|$. Moreover, the supremum in the definition of $ d_{H^{-1}}(G;E) $ is attained say by a function $f_G \in H^1(\pa E)$. By the condition \eqref{eq:zero-condition} for $\xi_{G,E}$ we may assume that the function $f_G $ has zero average on each component  $\Sigma_i$. Therefore we have by the Sobolev embedding 
\begin{equation}
\label{eq:supremum-d-G}
\|f_G\|_{L^p(\pa E)} \leq C_p \|\nabla_{\pa E} f_G\|_{L^2(\pa E)} \leq C_p 
\end{equation}
for all $p>2$.
By the minimality of $F$ it holds
\begin{equation}
\label{eq:from-minimality}
P(F) - P(G) \leq \Lambda_2 ( d_{H^{-1}}(G;E)- d_{H^{-1}}(F;E) )
\end{equation}
and
\begin{equation}
\label{eq:calculate-distance-diff}
	\begin{split}
	 d_{H^{-1}}&(G;E)- d_{H^{-1}}(F;E) 	\\
		&\leq \int_{\R^3} f_G(\pi_{\pa E}(x)) (\chi_{G}(x)-\chi_E(x))dx- \int_{\R^3} f_G(\pi_{\pa E}(x)) (\chi_{F}(x)-\chi_E(x))dx\\
		&=  \int_{\R^3}f_G(\pi_{\pa E}(x)) (\chi_{G}(x)-\chi_{F}(x))dx \\
		&\leq  \| f_G\circ  \pi_{\pa E}  \|_{L^p(\mathcal{N}_{\delta}(\pa E))}  \|\chi_{G}-\chi_{F}  \|_{L^{\frac{p}{p-1}}(\R^3)}.
	\end{split}
\end{equation}
The inequality \eqref{eq:supremum-d-G} implies  
$ \| f_G\circ  \pi_{\pa E}  \|_{L^p(\mathcal{N}_{\delta}(\pa E))}  \leq C$ and therefore we deduce from  \eqref{eq:from-minimality} and \eqref{eq:calculate-distance-diff}  that 
\[
P(F) \leq P(G) +  C \Lambda_2 |G \Delta F|^{1 - \frac{1}{p}}. 
\]
Hence, the minimizer $F$ is a constrained  $(\Lambda, \alpha)$-minimizer in the sense of Definition~\ref{def:constr-lambda-mini} for $\alpha =\frac{1}{p}$.  By  Lemma \ref{lem:constraint-omega-mini}, $F$ is also 
$(\Lambda_1, \alpha)$-minimizer for $\Lambda_1$ depending on $r_0, R_0, p$. 

By choosing $\delta $ small enough in \eqref{def:pena-functional},  we may deduce by Lemma \ref{lem:omega-mini}  that $F$ is $C^{1}$-regular and 
there is $\psi \in C^{1}(\pa E)$ with   $\|\psi\|_{C^1(\pa E)} \leq \delta_1$, where $\delta_1$ is as in \eqref{psi-delta1}, such that 
\[
\pa F = \{ x + \psi(x)\nu_E(x) : x \in \pa E\}.
\]
We then have from \eqref{eq:geometric-2}, which we proved above, that 
\[
P(E) \leq P(F)+ \Lambda_2 d_{H^{-1}}(F;E). 
\]
This proves that $E$ is a minimizer of \eqref{def:pena-functional}, which implies the claim. 
 
\end{proof}

\section{Regularity estimates for the flow}
\subsection{Definition of the flat flow solution}
We begin by giving   the construction of the flat flow solution for the surface diffusion flow, which is a slight modification of the one proposed by  Cahn-Taylor in \cite{CaTa} by introducing an additional constraint.   We first fix our constraint parameter $\delta>0$ and then  the time step $h>0$, and  given a bounded set  $E \subset \R^3$ of finite perimeter, which coincides with its Lebesgue representative, we consider the minimization problem 
\beq
\label{def:min-prob}
\min   \Big{\{} P(F)+\frac{d_{H^{-1}}(F;E)^2}{2h} : \quad   F \Delta E \subset  \overline{\mathcal{N}_{\delta}}(\pa^* E) \Big{\}}  
\eeq
and note that the minimizer exists. The constraint in \eqref{def:min-prob}, which requires the sets to be in a neighborhood of $E$, simplifies the forthcoming analysis, and as we will see, it does not affect the actual 
construction if the set $E$ is regular enough. Indeed, we will see that the minimizer $F$ of \eqref{def:min-prob} is  also regular and satisfies $ F \Delta E \subset \mathcal{N}_{h^{\gamma}}(\pa^* E)$ for $\gamma>0$ 
when $h$ is small enough. In other words,  we need the constraint in  \eqref{def:min-prob} to overcome the lack of coerciveness with respect to $d_{H^{-1}}$.

Let $E_0 \subset \R^3$  be a bounded set of finite perimeter  which coincides with its Lebesgue representative. 
We construct  the discrete-in-time evolution $\{E_k^{h,\delta}\}_{k\in \N}$ recursively setting $E(0)=E_0^{h,\delta}$ and, assuming that $E_k^{h,\delta}$ is defined, we set $E_{k+1}^{h,\delta}$ to be a minimizer of \eqref{def:min-prob} with $E= E_k^{h,\delta}$.  We define the \emph{ approximate flat flow}  $\{E^{h,\delta}(t)\}_{t \geq 0}$ by setting
\[
E^{h,\delta}(t) = E_k^{h,\delta} \qquad \text{for }\, t \in [kh, (k+1) h).
\]
We define a flat flow solution $\{E^\delta(t)\}_{t \geq 0}$ of the surface diffusion as any cluster point when we let  $h \to 0$.

\subsection{The first regularity estimates} 

Let us  study the minimization problem \eqref{def:min-prob}.  Throughout the rest of this section we will assume that the set $E \subset \R^3$ is $C^5$-regular  and satisfies the assumptions of  Proposition \ref{prop:geometric}, i.e., 
\begin{equation}
\label{eq:assumptions-E}
E \, \text{ satisfies the UBC with radius } \, r_0, \quad E \subset B_{R_0} \quad \text{and} \quad \|\Delta_{\pa E }H_E\|_{L^2(\pa E)} \leq K_0, 
\end{equation}
for some $K_0>1$. Note that the above assumptions imply that the set $E$  is uniformly $C^{2,\gamma}$-regular for any $\gamma \in (0,1)$.  The first  consequence of Proposition  \ref{prop:geometric} is the following. 

\begin{proposition}
\label{prop:apirori-estimates}
Assume $E \subset \R^3$ is $C^5$-regular,  satisfies the assumptions \eqref{eq:assumptions-E}, and let $F\subset \R^3$ be a minimizer of \eqref{def:min-prob} with $\delta < \delta_1$, where $\delta_1 $ is from Proposition \ref{prop:geometric}.  Then there is $C_0$, depending on $r_0,R_0,K_0$, such that 
\begin{equation}
\label{eq:apriori-linear}
d_{H^{-1}}(F;E) \leq C_0 h
\end{equation} 
for  $h \leq h_0$, where $h_0$ depends on  $r_0,R_0,K_0$ and on $\delta$.  Moreover, $F$ is a $(\Lambda, \alpha)$-minimizer of the perimeter for every $\alpha \in (0, \tfrac13)$, where $\Lambda$ depends on  $r_0,R_0,K_0$ and $\alpha$. 

In addition, if $\delta < \delta_2$, where $\delta_2 =  \min\{\delta_0,\delta_1\}$ and  $\delta_0$  is from Lemma~\ref{lem:omega-mini} with $\eps = \delta_1$,   then   $F $ is $C^5$-regular,  $\|H_F\|_{H^1(\pa F)} \leq C_0$,  there is $\psi \in C^1(\pa E)$  such that $\pa F = \{ x + \psi(x)\nu_E(x) : x \in \pa E \}$ and we have the estimates 
\begin{equation}
\label{eq:apriori-regular}
\|\psi\|_{H^1(\pa E)} \leq C_0 \sqrt{h} \quad \text{and} \quad \|\psi\|_{C^1(\pa E)} \leq h^{\frac{1}{9}},
\end{equation}
 when $h \leq h_0$. Moreover, for any $\eps>0$ there is $C_\eps$ such that $\|\psi\|_{L^\infty(\pa E)} \leq C_\eps h^{\frac{1}{2} -\eps}$. 
\end{proposition}

\begin{proof}
The estimate  \eqref{eq:apriori-linear} follows immediately by using the minimality of $F$ against $E$, and from \eqref{eq:geometric-1} stated in Proposition \ref{prop:geometric}, 
\[
P(F)+\frac{d_{H^{-1}}(F;E)^2}{2h} \leq P(E) \leq P(F) + \Lambda_2 d_{H^{-1}}(F;E). 
\]

Next we prove that the minimizer $F$ is  a $(\Lambda,\alpha)$-minimizer of the perimeter for any $\alpha \in (0,\tfrac13)$. Arguing as in the proof of Proposition \ref{prop:geometric},  we may 
reduce to  considering sets  $G \subset \overline{\mathcal{N}_{\delta}}(E) $ such that $G \Delta F \subset \overline{\mathcal{N}_{\delta}}(\Sigma_i)$ where $\Sigma_i$ is a component of $\pa E$. For these sets we have by Lemma   \ref{lem:distance-finite}  that $ d_{H^{-1}}(G;E)  <\infty$ if  $|G| = |F|$. Moreover, the supremum in the definition of $ d_{H^{-1}}(G;E)$ is attained say by a function $f_G \in H^1(\pa E)$. By the condition \eqref{eq:zero-condition} for $\xi_{G,E}$ we may assume that the function $f_G $ has zero average on every component  $\Sigma_i$. Therefore we have by the Sobolev embedding  $\|f_G\|_{L^p(\pa E)} \leq C \|\nabla_{\pa E} f_G\|_{L^2(\pa E)} \leq C$ for $p = \frac{1}{\alpha}$. We will show that 
\begin{equation}
\label{eq:apirori-1}
P(F) \leq P(G) + \Lambda |G \Delta F|^{1-\alpha}.
\end{equation}
We remark that by Lemma \ref{lem:constraint-omega-mini}, the estimate \eqref{eq:apirori-1} then implies that $F$ is  a $(\Lambda_1, \alpha)$-minimizer of the perimeter according to the Definition \ref{def:lambda-mini}. We divide the proof of \eqref{eq:apirori-1} in two cases. 

Let us first  assume  that $d_{H^{-1}}(G;E) \leq 2C_0 h$, where $C_0$ is the constant from \eqref{eq:apriori-linear}. By the minimality of $F$ we have
\[
P(F)+\frac{d_{H^{-1}}(F;E)^2}{2h}\leq P(G)+\frac{d_{H^{-1}}(G;E)^2}{2h}.
\]
If $d_{H^{-1}}(G;E) \leq d_{H^{-1}}(F;E)$, then \eqref{eq:apirori-1} follows immediately. Let us then  assume  $d_{H^{-1}}(G;E)> d_{H^{-1}}(F;E)$.  By   \eqref{eq:apriori-linear}  it holds 
\begin{equation}
\label{eq:apirori-2}
\begin{split}
P(F) - P(G) &\leq \frac{1}{2h}\left(d_{H^{-1}}(G;E) + d_{H^{-1}}(F;E)\right) \left(d_{H^{-1}}(G;E) - d_{H^{-1}}(F;E)\right)\\
&\leq 2C_0 \left(d_{H^{-1}}(G;E) - d_{H^{-1}}(F;E)\right). 
\end{split}
\end{equation}
Using \eqref{eq:calculate-distance-diff}  we deduce 
\begin{equation}
\label{eq:apirori-3}
d_{H^{-1}}(G;E) - d_{H^{-1}}(F;E) \leq  \| f_G\circ  \pi_{\pa E}  \|_{L^p(\mathcal{N}_{\delta}(\pa E))}  \|\chi_{G}-\chi_{F}  \|_{L^{\frac{p}{p-1}}(\R^3)} \leq C |G \Delta F |^{1- \alpha}, 
\end{equation}
where the last inequality follows from $\|f_G\|_{L^p(\pa E)} \leq C$ and $p = \frac{1}{\alpha}$.  This together with \eqref{eq:apirori-2} implies \eqref{eq:apirori-1}. 

Let us then assume that  $d_{H^{-1}}(G;E) > 2C_0 h$. Note that by  \eqref{eq:apriori-linear}   it holds
\begin{equation}
\label{eq:apirori-4}
d_{H^{-1}}(F;E) < \frac12 d_{H^{-1}}(G;E).  
\end{equation}
Using first  the minimality of $F$ against $E$,  then   \eqref{eq:geometric-1} with $F$   replaced by $G$, and finally \eqref{eq:apirori-4},   we obtain
\[
\begin{split}
P(F) \leq P(E) &\leq P(G) +\Lambda_2 d_{H^{-1}}(G;E)\\
&\leq  P(G) + 2 \Lambda_2 \big(d_{H^{-1}}(G;E) - d_{H^{-1}}(F;E)\big).
\end{split}
 \]
The claim \eqref{eq:apirori-1} then follows from \eqref{eq:apirori-3}.

We may use Lemma \ref{lem:omega-mini} to deduce that when $\delta$ is small enough,  there  is $\psi \in C^{1,\frac13}(\pa E)$  with $\|\psi\|_{C^{1,\frac13}(\pa E)} \leq \delta_1$, where $\delta_1$ is from Proposition \ref{prop:geometric}, such that $\pa F = \{ x + \psi(x)\nu_E(x) : x \in \pa E \}$. Then \eqref{eq:geometric-2}, the minimality of $F$ and  \eqref{eq:apriori-linear} imply
\[
\frac14\| \nabla_{\pa E} \psi \|_{L^2(\pa E)}^2 +P(E) \leq P(F)+  \Lambda_2 d_{H^{-1}}(F;E) \leq P(E) +  Ch. 
\]
This yields
\[
\|\nabla_{\pa E} \psi \|_{L^2(\pa E)} \leq C \sqrt{h}.
\]
The estimate  $\|\psi\|_{H^1(\pa E)} \leq C \sqrt{h}$  then follows from \eqref{eq:geometric-44} and \eqref{eq:geometric-47}. In particular, we deduce by the Sobolev embedding and  the $C^1$-bound on $\psi$  that $\|\psi\|_{L^\infty(\pa E)} \leq C_\eps h^{\frac{1}{2} -\eps}$ for any  $\eps>0$ and therefore the minimizer $F$ does not touch the constraint $ \mathcal{N}_{\delta}(\pa E)$ when $h$ is small enough. Hence, we may  calculate the Euler-Lagrange equation for $F$ which states that 
\begin{equation}
\label{eq:Euler-Lagrange-1}
H_F + \frac{d_{H^{-1}}(F;E)}{h} f\circ \pi_{\pa E} = \lambda \qquad \text{on } \, \pa F,
\end{equation}
where $f$ is the function which attains the supremum in $d_{H^{-1}}(F;E)$. Note that since $E$ is $C^5$-regular and since $f$ solves \eqref{eq:akka-meno-uno}, we may bootstrap the regularity of $F$ to $C^{5,\gamma}$ for any $\gamma\in (0,1)$.   We remark that the $C^5$-regularity depends on $h$. 

The estimate $\|H_F\|_{H^1(\pa F)} \leq C$ then follows from \eqref{eq:Euler-Lagrange-1}, from $d_{H^{-1}}(F;E) \leq C_0h$,  and from the fact that $\|\nabla_{\pa E} f\|_{L^2(\pa E)} \leq 1$ which in turn yields $\|\nabla_{\pa F} (f\circ \pi_{\pa E}) \|_{L^2(\pa F)} \leq C$ by \eqref{eq:change-variables-function}.  Finally the estimate  $\|\psi\|_{C^1(\pa E)} \leq  h^{\frac{1}{9}}$ follows from   $\|\psi\|_{L^\infty(\pa E)} \leq C_\eps h^{\frac{1}{2} -\eps}$,  $\|\psi\|_{C^{1,\frac13}(\pa E)} \leq \delta_1$ and the  interpolation inequality \eqref{eq:Holder-interpolate}.  
\end{proof}

As we already mentioned, an important consequence of Proposition \ref{prop:apirori-estimates} is that a minimizer of  \eqref{def:min-prob} does not touch the constraint ${\mathcal{N}_{\delta}}(\pa^* E) $ when $h \leq h_0$. We may then calculate the Euler-Lagrange equation for  the minimizer $F$, and  obtain by  \eqref{eq:akka-meno-uno} and  \eqref{eq:Euler-Lagrange-1} that 
\begin{equation}\label{eq:Euler-Lagrange-3}
\begin{cases}
&H_F + \frac{d_{H^{-1}}(F;E)}{h} f\circ \pi_{\pa E} = \lambda \qquad \text{on } \, \pa F \\
&-\Delta_{\pa E} f =\frac{\xi_{F,E}}{d_{H^{-1}}(F;E)} \qquad \text{on } \pa E, 
\end{cases}
\end{equation}
where $\xi_{F,E}$ is defined in \eqref{def:normal-velo}. 
We remark that if $E$ is $C^{5,\gamma}$-regular for some $\gamma \in (0,1)$, then by the elliptic regularity theory we may conclude that $F$ is $C^{6,\gamma}$-regular. This means that the minimizing movement scheme has a regularizing effect. 

We may  use Lemma \ref{lem:expand-mean-curv} to combine the two above equations into the following 
\begin{equation}
\label{eq:Euler-Lagrange-2}
\frac{1}{h}\left( \psi + \frac{H_E}{2} \psi^2  + \frac{K_E}{3}\psi^3\right) = \Delta_{\pa E} \left( - \Delta_{\pa E}\psi  + H_E  \right) +  \Delta_{\pa E} R_0 \qquad \text{on } \, \pa E,
\end{equation}
where the error term $R_0$  is of the form \eqref{eq:error-R-0}. The estimates in  Proposition \ref{prop:apirori-estimates} do not guarantee that we keep  the uniform ball condition, because the bound $\|H_F\|_{H^1(\pa F)} \leq C$ does not imply a bound for $\|B_F\|_{L^\infty(\pa F)}$. Therefore we need to upgrade the regularity estimates from Proposition \ref{prop:apirori-estimates}, which we do by using the equation \eqref{eq:Euler-Lagrange-2}.

\begin{lemma}
\label{lem:upgrade-regularity}
Assume $E \subset \R^3$ is $C^5$-regular, satisfies the assumptions \eqref{eq:assumptions-E} and $\|\Delta_{\pa E} H_E\|_{H^1(\pa E)}$ $ \leq  K_0 h^{-\frac14}$. Let  $F\subset \R^3$ be a minimizer of \eqref{def:min-prob} with $\delta < \delta_2$ where $\delta_2$ is from Proposition \ref {prop:apirori-estimates}.  Then for the heightfunction in \eqref{eq:apriori-regular} it holds
\begin{equation}
\label{eq:upgrade-regularity}
\|\psi \|_{L^2(\pa E)} \leq C_1h  \quad \text{and} \quad  \|\psi \|_{H^4(\Sigma)} \leq C_1 
\end{equation}
for all $h \leq h_0$, where $h_0$ is from Proposition \ref {prop:apirori-estimates} and $\Sigma = \pa E$. The constant $C_1$ depends on $r_0, R_0$ and  $K_0$.
\end{lemma}
\begin{proof}
We begin by noticing that the assumptions on $E$ and Lemma \ref{lem:hilbert-norm} imply 
\begin{equation}
\label{eq:upgr-reg-0}
  \|B_\Sigma \|_{L^\infty(\Sigma)} \leq \hat C, \quad  \|B_\Sigma \|_{H^2(\Sigma)} \leq \hat C \, K_0 \quad \text{and} \quad  \|B_\Sigma\|_{H^3(\Sigma)} \leq  \hat C\, K_0  h^{-\frac14},
\end{equation}
for a constant $\hat C$, which depends on $r_0$ and $R_0$.

We  obtain the  estimates by multiplying the Euler-Lagrange equation \eqref{eq:Euler-Lagrange-2} with $\Delta_{\pa E}^2 \psi $. Integrating by parts and by  re-organizing the terms we deduce 
\[
\begin{split}
\frac{1}{h} \int_{\pa E } |\Delta_{\pa E} \psi|^2 \, d \H^2 +  \int_{\pa E } |\Delta_{\pa E}^2 \psi|^2 \, d \H^2 &=  \int_{\pa E } (\Delta_{\pa E}^2 \psi) (\Delta_{\pa E} H_E + \Delta_{\pa E} R_0 ) \, d \H^2\\
&- \frac{1}{h} \int_{\pa E } (\Delta_{\pa E} \psi) \Delta_{\pa E} ( \frac{H_E}{2} \psi^2  + \frac{K_E}{3}\psi^3 )  \, d \H^2.
\end{split}
\]
By Cauchy-Schwarz inequality this implies 
\begin{equation}
\label{eq:upgr-reg-1}
\begin{split}
 \frac{1}{h}  \|\Delta_{\pa E} \psi\|_{L^2}^2 +  \|\Delta_{\pa E}^2 \psi\|_{L^2}^2 \leq  &\|\Delta_{\pa E} H_E\|_{L^2}^2  +  \|\Delta_{\pa E} R_0\|_{L^2}^2  \\
&+ \frac{1}{h}  \| \Delta_{\pa E}  (\frac{H_E}{2} \psi^2)  \|_{L^2}^2 +\frac{1}{h}  \| \Delta_{\pa E}( \frac{K_E}{3}\psi^3 )\|_{L^2}^2 .  
\end{split}
\end{equation}
We need  to estimate all four terms on  the right-hand-side in \eqref{eq:upgr-reg-1}. 

By  the assumptions  \eqref{eq:assumptions-E} it holds 
\begin{equation}
\label{eq:upgr-reg-2}
\|\Delta_{\pa E} H_E\|_{L^2}^2  \leq K_0^2. 
\end{equation}

The term $ \|\Delta_{\pa E} R_0\|_{L^2}^2$ is the most difficult to deal with.   In order to estimate it, we recall the form of $R_0$ in \eqref{eq:error-R-0}. Therefore we have by the Leibniz rule
\[
\begin{split}
\|\Delta_{\pa E} R_0\|_{L^2(\Sigma)} \leq \,   &C\sum_{j+k = 2} \|  |\bar \nabla^{j}( A_1(\psi B_\Sigma ,\bar \nabla \psi ) | | \bar \nabla^{2+k} \psi | \|_{L^2(\Sigma)} \\
&+ C\sum_{j+k = 2}   \| |\bar \nabla^{j}  A_2(\psi B_\Sigma,\bar \nabla \psi )| |  \bar \nabla^{1+k} (\psi B_\Sigma)|\|_{L^2(\Sigma)}\\
&+ \| a_0(\psi,\bar \nabla \psi,  B_\Sigma)\|_{H^2(\Sigma)}.
\end{split}
\]  
For $j,k$ with $j+k = 2$  we use Lemma  \ref{lem:Leibniz} with $T_1 = A_1(\psi B_\Sigma ,\bar \nabla \psi )$  and $T_2 = \bar \nabla \psi$ to estimate 
\begin{equation}
\label{eq:upgr-reg-3}
\begin{split}
 \|  |\bar \nabla^{j}( A_1(\psi B_\Sigma ,\bar \nabla \psi ) | &| \bar \nabla^{2+k} \psi | \|_{L^2} \\
&\leq C \|A_1(\psi B_\Sigma ,\bar \nabla \psi )\|_{L^\infty} \|\psi \|_{H^4} + C\|\psi\|_{C^1} \|A_1(\psi B_\Sigma ,\bar \nabla \psi )\|_{H^3}
\end{split}
\end{equation}
and similarly with  $T_1 = A_2(\psi B_\Sigma ,\bar \nabla \psi )$  and $T_2 =\psi B_\Sigma$ 
\begin{equation}
\label{eq:upgr-reg-4}
\begin{split}
 \|  |\bar \nabla^{j}( A_2(\psi B_\Sigma ,\bar \nabla \psi )|  &|  \bar \nabla^{1+k} (\psi B_\Sigma)|  \|_{L^2}\\
&\leq C \| A_2(\psi B_\Sigma ,\bar \nabla \psi )\|_{L^\infty} \|\psi B_\Sigma \|_{H^3} + C\|\psi B_\Sigma \|_{L^\infty} \| A_2(\psi B_\Sigma ,\bar \nabla \psi )\|_{H^3}.
\end{split}
\end{equation}
Since $A_i$ is smooth and $A_i(0,0) = 0$, it holds 
\begin{equation}
\label{eq:upgr-reg-6}
 \|A_i(\psi B_\Sigma ,\bar \nabla \psi )\|_{L^\infty} \leq C \|\psi\|_{C^1}. 
\end{equation} 
Moreover, again by the smoothness of $A_i$ we have by the chain rule that   
\[
\begin{split}
| \bar \nabla^3 ( A_i(\psi B_\Sigma ,\bar \nabla \psi )) |  \leq C \sum_{|\alpha| \leq 3}&(1+  |\bar \nabla^{\alpha_1} (\psi B_\Sigma)|)(1+   |\bar \nabla^{\alpha_2} (\psi B_\Sigma)|)(1+|\bar \nabla^{\alpha_3} (\psi B_\Sigma)| )\cdots \\
 &\cdots(1+  |\bar \nabla^{1+\alpha_4} \psi|)(1+   | \bar \nabla^{1+\alpha_5} \psi|)(1+ | \bar \nabla^{1+\alpha_6} \psi| ).
\end{split}
\]
 Therefore Lemma  \ref{lem:Leibniz}  for $T_i = \psi B_\Sigma$  for $i=1,2,3$ and   $T_i =\bar \nabla \psi $ for $i =4,5,6$  implies 
\begin{equation}
\label{eq:upgr-reg-7}
\begin{split}
 \|A_i(\psi B_\Sigma ,\bar \nabla \psi )\|_{H^3} &\leq C\big( (1+ \|\psi B_\Sigma \|_{L^\infty})(1+   \|\psi \|_{H^4})  + (1+ \|\psi  \|_{C^1})(1+\|\psi B_\Sigma \|_{H^3})\big)\\
&\leq C (1+\|\psi \|_{H^4}  + \|\psi B_\Sigma \|_{H^3}).
\end{split}
\end{equation}
Finally, we estimate the third term similarly by the chain rule, by the regularity of $a_0$,  \eqref{eq:upgr-reg-0},   Lemma  \ref{lem:Leibniz}  and by interpolation  inequality in Proposition \ref{prop:interpolation} as
\begin{equation}
\label{eq:upgr-reg-5}
\begin{split}
 \|a_0(\psi,\bar \nabla \psi,  &B_\Sigma)\|_{H^2(\Sigma)} \\
&\leq   C\sum_{|\alpha| \leq 2} \| (1+ |\bar \nabla^{1+\alpha_1}\psi  |)  (1+|\bar \nabla^{1+\alpha_2}\psi  |)(1+ | \bar \nabla^{\alpha_3} B_\Sigma | )(1+ | \bar \nabla^{\alpha_4} B_\Sigma |)\|_{L^2(\Sigma)} \\
&\leq C(1+ \| \psi\|_{H^3} + \|B_\Sigma\|_{H^2} ) \\
&\leq C(1+\| \psi\|_{H^3}) \leq \eps \| \psi\|_{H^4} + C_\eps  \| \psi\|_{H^2} +C.
\end{split}
\end{equation}

  We  may thus estimate  the term $\|\Delta_{\pa E} R_0\|_{L^2(\Sigma)}$ by    \eqref{eq:upgr-reg-3}, \eqref{eq:upgr-reg-4}, \eqref{eq:upgr-reg-6}, \eqref{eq:upgr-reg-7} and  \eqref{eq:upgr-reg-5} as
\[
\|\Delta_{\pa E} R_0\|_{L^2(\Sigma)} \leq C  \|\psi\|_{C^1} \big( 1 +\|\psi \|_{H^4} +   \|\psi B_\Sigma \|_{H^3} \big) +\eps \| \psi\|_{H^4} + C_\eps  \| \psi\|_{H^2} +C.
\]
Using  Lemma  \ref{lem:Leibniz} once again we obtain $ \|\psi B_\Sigma \|_{H^3}  \leq C \|\psi\|_{L^\infty }    \|B_\Sigma \|_{H^3}  +  C\|\psi \|_{H^3}$. We have by Proposition \ref{prop:apirori-estimates}  that 
$ \|\psi\|_{L^\infty }  \leq C_\eps h^{\frac12-\eps}$ for any $\eps>0$ and by  \eqref{eq:upgr-reg-0}  we have $\|B_\Sigma\|_{H^3(\Sigma)} \leq C_1 h^{-\frac14}$, which yield 
\begin{equation}
\label{eq:upgr-reg-88}
\|\psi\|_{L^\infty }    \|B_\Sigma \|_{H^3}  \leq C_\eps h^{\frac14-\eps}  \leq 1,
\end{equation}
when $\eps <\frac14$ and  $h$ is small.  Therefore we may bound  $\|\Delta_{\pa E} R_0\|_{L^2(\Sigma)}$ by the previous estimates as
\[
\|\Delta_{\pa E} R_0\|_{L^2(\pa E)} \leq C \|\psi\|_{C^1}  \|\psi \|_{H^4}  +\eps \| \psi\|_{H^4} + C_\eps  \| \psi\|_{H^2} +C.
\]
Finally, we use  Lemma \ref{lem:hilbert-norm} with $k=1,2$ and \eqref{eq:upgr-reg-0}  to deduce that 
\begin{equation}
\label{eq:upgr-reg-8-end}
 \|\psi \|_{H^2(\Sigma)} \leq C (1 + \|\Delta_{\pa E} \psi \|_{L^2(\pa E)} ) \quad \text{and} \quad  \|\psi \|_{H^4(\Sigma)} \leq C (1 + \|\Delta_{\pa E}^2\psi \|_{L^2(\pa E)} ).
\end{equation}
Hence, we obtain   by  \eqref{eq:apriori-regular}   and choosing $\eps$ small enough 
\begin{equation}
\label{eq:upgr-reg-8}
\begin{split}
\|\Delta_{\pa E} R_0\|_{L^2(\pa E)} \leq \frac12 \|\Delta_{\pa E}^2\psi \|_{L^2(\pa E)} +C \|\Delta_{\pa E} \psi \|_{L^2(\pa E)}  +C.
\end{split}
\end{equation}

We are left with the two last terms in \eqref{eq:upgr-reg-1}. We estimate only the second last term  $\frac{1}{h}  \| \Delta_{\pa E}  (\frac{H_E}{2} \psi^2)  \|_{L^2}^2$ as the final term is treated with the same argument. We first use Lemma~\ref{lem:Leibniz},  \eqref{eq:upgr-reg-0}, the Sobolev embedding, interpolation inequalities and \eqref{eq:hilbert-norm}  to estimate
\begin{equation}
\label{eq:upgr-reg-88}
\begin{split}
\| \Delta_{\pa E}  (\frac{H_E}{2} \psi^2)  \|_{L^2} &\leq C \|\psi\|_{L^\infty}^2 \|H_\Sigma\|_{H^2 } + C \|\psi\|_{L^\infty} \|H_\Sigma \|_{L^\infty} \|\psi\|_{H^2 } +C \|H_\Sigma \|_{L^\infty}\|\nabla_{\pa E}\psi\|^2_{L^4}\\
&\leq C  \|\psi\|_{L^\infty} \|\psi\|_{H^2 } \\
&\leq C  \|\psi\|_{L^\infty}( \| \Delta_{\pa E}  \psi\|_{L^2} + C  \|\psi\|_{L^2}).
\end{split}
\end{equation}
Proposition \ref{prop:apirori-estimates} implies $\|\psi \|_{L^\infty} \leq C_\eps h^{\frac12 -\eps}$. We choose $\eps= \frac14$ and deduce
\begin{equation}
\label{eq:upgr-reg-9}
\frac{1}{h}\| \Delta_{\pa E}  (\frac{H_E}{2} \psi^2)  \|_{L^2}^2 \leq \frac{C}{\sqrt{h}}  \| \Delta_{\pa E}  \psi\|_{L^2}^2 + C. 
\end{equation}
The same argument also yields 
\begin{equation}
\label{eq:upgr-reg-10}
\frac{1}{h}  \| \Delta_{\pa E}( \frac{K_E}{3}\psi^3 )\|_{L^2}^2 \leq \frac{C}{\sqrt{h}}  \| \Delta_{\pa E}  \psi\|_{L^2}^2 + C.
\end{equation}

We obtain  from \eqref{eq:upgr-reg-1}, \eqref{eq:upgr-reg-2}, \eqref{eq:upgr-reg-8},  \eqref{eq:upgr-reg-9} and  \eqref{eq:upgr-reg-10} that 
\[
\frac{1}{h}  \|\Delta_{\pa E} \psi\|_{L^2}^2 +  \|\Delta_{\pa E}^2 \psi\|_{L^2}^2 \leq  C  +  \frac12 \|\Delta_{\pa E}^2\psi \|_{L^2}^2 + \frac{C}{\sqrt{h}}  \| \Delta_{\pa E}  \psi\|_{L^2}^2
\]
when $h$ is small. This yields
\[
\frac{1}{h}  \|\Delta_{\pa E} \psi\|_{L^2}^2 +  \|\Delta_{\pa E}^2 \psi\|_{L^2}^2 \leq  C,
\]
when $h$ is small. Therefore we obtain $\|\psi\|_{H^4} \leq C$ using \eqref{eq:upgr-reg-8-end}. 

In order to prove $\|\psi\|_{L^2} \leq C h$, we remark that the previous argument, in fact, implies that the right-hand-side of \eqref{eq:Euler-Lagrange-2} is bounded in $L^2(\pa E)$. Therefore we deduce that the left-hand-side of  \eqref{eq:Euler-Lagrange-2} is also bounded in $L^2(\pa E)$, i.e.,
\[
\frac{1}{h} \|\xi_{F,E}\|_{L^2(\pa E)} =  \frac{1}{h} \|\psi + \frac{H_E}{2} \psi^2  + \frac{K_E}{3}\psi^3\|_{L^2(\pa E)} \leq C. 
\]
The estimate $\|\psi\|_{L^2} \leq C h$ then follows from \eqref{eq:geometric-44}.

\end{proof}

The previous lemma implies regularity for the heighfunction. We also need an opposite result which  states that regularity for the heighfunction implies regularity for the minimizing set $F$. 
\begin{lemma}
\label{lem:geometric-regularity}
Assume $E \subset \R^3$ is $C^5$-regular, satisfies the assumptions \eqref{eq:assumptions-E} and  $\|\Delta_{\pa E} H_E\|_{H^1(\pa E)} \leq K_0 h^{-\frac14}$ for $K_0\geq 1$. Let   $F\subset \R^3$ be a minimizer of \eqref{def:min-prob} with $\delta < \delta_2$, where $\delta_2$ is from Proposition \ref {prop:apirori-estimates}. Assume that  the heightfunction in \eqref{eq:apriori-regular} satisfies
\[
\|\psi\|_{L^2(\pa E)} \leq L_0 h \qquad \text{and} \qquad  \|\Delta_{\pa E} \psi\|_{L^2(\pa E)} \leq L_0 \sqrt{h}. 
\]
Then there is a constant $C_2$, depending on $r_0, R_0$ and $L_0$, such that  
\[
\|\Delta_{\pa F} H_F \|_{L^2(\pa F)} \leq C_2(1 + \sqrt{K_0})  \qquad \text{and} \qquad  \|\nabla_{\pa F} \Delta_{\pa F} H_F \|_{L^2(\pa F)}  \leq C_2 h^{-\frac14}
\]
for $h \leq h_0$, where $h_0$ depends on $r_0, R_0, K_0$ and $\delta$. Moreover, $F$ satisfies the UBC with radius $r_0/2$.
\end{lemma}

\begin{proof}
Below we denote by $C$ a generic constant that depends on $r_0, R_0, L_0$. We need only to be careful to trace the dependence on $K_0$. 

 Note first that by \eqref{eq:hilbert-norm} and by interpolation from Proposition \ref{prop:interpolation}  it holds 
\begin{equation}\label{eq:geometric-regularity-1}
\|\psi\|_{H^1(\pa E)} \leq C h^{\frac34} . 
 \end{equation}
Moreover,  by Lemma \ref{lem:upgrade-regularity} it holds  $\|\psi \|_{H^4(\Sigma)} \leq C_1$ for $C_1 = C_1(r_0,R_0,K_0)$, which together with \eqref{eq:geometric-regularity-1}, with Sobolev embedding and interpolation inequality \eqref{eq:Holder-interpolate} yields that given $0<\alpha<1$ there exists $\gamma>0$ such that $\|\psi \|_{C^{2,\alpha}(\Sigma)} \leq C_1h^{\gamma}$ where $C_1 = C_1(r_0,R_0,K_0,\alpha)$. From here we infer, see e.g.  \cite[Theorem 2.6]{Dal} and \cite[Remark 3.4.6]{Antonia}, that $F$ satisfies the UBC with radius $r_0/2$ when $h$ is small. In particular, 
\begin{equation}\label{eq:geometric-regularity-33}
\|B_F\|_{L^\infty(\pa F)} \leq \frac{2}{r_0}. 
 \end{equation}

We will obtain the estimates again from the Euler-Lagrange equations \eqref{eq:Euler-Lagrange-3}.  We define $\tilde f = \frac{d_{H^{-1}}(F;E)}{h} f$ and write \eqref{eq:Euler-Lagrange-3} as
\begin{equation}\label{eq:Euler-Lagrange-4}
\begin{cases}
&H_F +  (\tilde f\circ \pi_{\pa E}) = \lambda \qquad \text{on } \, \pa F \\
&-\Delta_{\pa E} \tilde f =\frac{\xi_{F,E}}{h} \qquad \text{on } \pa E, 
\end{cases}
\end{equation}
where $\xi_{F,E}$ is defined in \eqref{def:normal-velo}.

By \eqref{eq:geometric-44}, by the assumption $\|\psi\|_{L^2(\pa E)} \leq L_0 h$  and by \eqref{eq:geometric-regularity-1}  it holds 
\begin{equation}\label{eq:geometric-regularity-2}
\|\xi_{F,E}\|_{L^2(\pa E)} \leq C h  \quad \text{and} \quad  \|\xi_{F,E}\|_{H^1(\pa E)} \leq C h^{\frac34}. 
\end{equation}
We have by \eqref{eq:hilbert-norm},   by Lemma  \ref{lem:hilbert-norm}, by  \eqref{eq:upgr-reg-0} and by the second equation in  \eqref{eq:Euler-Lagrange-4}   that  
\begin{equation}\label{eq:geometric-regularity-3}
\begin{split}
&\|\tilde f \|_{H^2(\Sigma)} \leq C(1 + \|\Delta_{\pa E} \tilde f\|_{L^2(\pa E)})  \leq C(1+ \frac{\|\xi_{F,E}\|_{L^2(\pa E)}}{h}) \leq C \qquad \text{and}\\
&\|\tilde f \|_{H^3(\Sigma)} \leq C(1 + \|B_{\Sigma}\|_{H^2(\Sigma)} + \| \Delta_{\pa E} \tilde f\|_{H^1(\pa E)})  \leq C(K_0+\frac{\|\xi_{F,E}\|_{H^1(\pa E)}}{h}) \leq C h^{-\frac14},
\end{split}
\end{equation}
when $h \leq h_{0}(r_0, R_0, K_0,\delta)$ is small. Above we also used the fact that $\tilde f$ has zero average  in every component of $\pa E$ and Sobolev embedding. 

Next we need to relate the derivative of  the function $\tilde f\circ \pi_{\pa E}$ on $\Gamma = \pa F$ to the derivatives of $\tilde f$ on $\Sigma = \pa E$. First, we have by \eqref{eq:change-variables-function} that 
\[
|\nabla_{\pa F} (\tilde f\circ \pi_{\pa E}) (x)|  \leq C |\nabla_{\pa E} \tilde f (\pi_{\pa E}(x))| 
\]
for all $x \in \pa F$. We write this in terms of covariant derivatives as
\begin{equation}\label{eq:geometric-regularity-4}
|\bar \nabla_{\Gamma} (\tilde f\circ \pi_{\pa E}) (x)| \leq C |\bar \nabla_{\Sigma} \tilde f (\pi_{\pa E}(x))|
\end{equation}
for all $x \in \Gamma$.  The higher order version of \eqref{eq:geometric-regularity-3} is technically more involved and to that aim, we proceed by using   \cite[Lemma 5.3]{JN} with $\Sigma = \pa E$ to deduce for $k = 2,3$  and for all $x \in \Gamma = \pa F$ that 
\begin{equation}\label{eq:geometric-regularity-5}
|\nabla^{k} (\tilde f\circ \pi_{\pa E}) (x)| \leq C \sum_{|\alpha| \leq k} ( 1+ |\bar \nabla_\Sigma^{\alpha_1} B_\Sigma(\pi_{\pa E}(x))|) \cdots (1+ |\bar \nabla_\Sigma^{\alpha_{k-1}} B_\Sigma(\pi_{\pa E}(x))|) |\bar  \nabla_\Sigma^{\alpha_{k}} \tilde f (\pi_{\pa E}(x))|,
\end{equation}
 where $\bar  \nabla_\Sigma^j$ denotes   $j$-th order covariant derivative on $\Sigma$ and  $\nabla^j$ denotes the $j$-th differential in $\R^3$. Even if it is not explicitly stated in \cite[Lemma 5.3]{JN},  it is clear that in \eqref{eq:geometric-regularity-5}  the last index is  never zero, i.e., $\alpha_k \geq 1$.

Let us first prove the estimate $\|\Delta_{\pa F} H_F \|_{L^2(\pa F)} \leq C_2(1 + \sqrt{K_0})$.  We use \eqref{eq:geometric-regularity-5} for $k=2$ and \eqref{eq:geometric-regularity-33} to estimate point-wise on $x \in \pa F$ (recall that in \eqref{eq:geometric-regularity-5}  $\alpha_2 \geq 1$)
\[
|\nabla^{2} (\tilde f\circ \pi_{\pa E}) (x)| \leq C |\bar  \nabla_\Sigma^{2} \tilde f (\pi_{\pa E}(x))| + C( 1+ |\bar \nabla_\Sigma B_\Sigma(\pi_{\pa E}(x))|)  |\bar  \nabla_\Sigma \tilde f (\pi_{\pa E}(x))|.
\]
By Cauchy-Schwarz inequality, by the Sobolev embedding and by \eqref{eq:geometric-regularity-3}  it holds 
\[
\begin{split}
\|\nabla^{2} (\tilde f\circ \pi_{\pa E})\|_{L^2(\Gamma)} &\leq C( \| \tilde f \|_{H^2(\Sigma)} + (1+ \|\bar \nabla_\Sigma  B_\Sigma\|_{L^4(\Sigma)})  \|\bar \nabla_\Sigma \tilde f \|_{L^4(\Sigma)}  )\\
&\leq C( \| \tilde f \|_{H^2(\Sigma)} + (1+ \|\bar \nabla_\Sigma  B_\Sigma\|_{L^4(\Sigma)}) \| \tilde f \|_{H^2(\Sigma)} )\\
&\leq C(1 +  \|\bar \nabla_\Sigma  B_\Sigma\|_{L^4(\Sigma)}).
\end{split}
\]
By the interpolation inequality in Proposition \ref{prop:interpolation},  by \eqref{eq:upgr-reg-0} and by the fact that $|B_\Sigma| \leq C$   it holds 
\[
 \|\bar \nabla_\Sigma  B_\Sigma\|_{L^4(\Sigma)} \leq C  \| B_\Sigma\|_{H^2(\Sigma)}^{\frac12} \| B_\Sigma\|_{L^\infty(\Sigma)}^{\frac12}\leq C \sqrt{K_0}. 
\]
Therefore we have 
\begin{equation}\label{eq:geometric-regularity-7}
\|\nabla^{2} (\tilde f\circ \pi_{\pa E})\|_{L^2(\Gamma)} \leq C(1  + \sqrt{K_0}).
\end{equation}

We may write the Laplacian of $ \tilde f\circ \pi_{\pa E}$  on $\pa F$  as
\begin{equation}\label{eq:geometric-laplacian}
\Delta_{\pa F} (\tilde f\circ \pi_{\pa E})  = \Delta_{\R^3} (\tilde f\circ \pi_{\pa E}) - \langle\nabla^2 (\tilde f\circ \pi_{\pa E}) \nu_F, \nu_F\rangle - H_F (\nabla (\tilde f\circ \pi_{\pa E}) \cdot  \nu_F). 
\end{equation}
Therefore we have by \eqref{eq:geometric-regularity-33}  and \eqref{eq:geometric-regularity-7}  that 
\begin{equation}\label{eq:geometric-regularity-77}
\|\Delta_{\pa F} (\tilde f\circ \pi_{\pa E}) \|_{L^2(\pa F)} \leq C(1  + \sqrt{K_0}). 
\end{equation}
The first  equation  in \eqref{eq:Euler-Lagrange-4} then  implies  $\| \Delta_{\pa F} H_F \|_{L^2(\pa F)} \leq C_2(1  + \sqrt{K_0})$ and the claim follows.

We need yet to prove  $\| \nabla_{\pa F}\Delta_{\pa F} H_F \|_{L^2(\pa F)} \leq C_2 h^{-\frac14} $.  We use  \eqref{eq:geometric-regularity-5} for $k=3$  and  \eqref{eq:geometric-regularity-33}  to estimate point-wise on $x \in \pa F$
\[
\begin{split}
|\nabla^{3} (\tilde f\circ \pi_{\pa E}) (x)| \leq &C |\bar  \nabla_\Sigma^{3} \tilde f (\pi_{\pa E}(x))| + C( 1+ |\bar \nabla_\Sigma B_\Sigma(\pi_{\pa E}(x))|)  |\bar  \nabla_\Sigma^2 \tilde f (\pi_{\pa E}(x))| \\
&+C( 1+ |\bar \nabla_\Sigma^2 B_\Sigma(\pi_{\pa E}(x))| + |\bar \nabla_\Sigma B_\Sigma(\pi_{\pa E}(x))|^2)  |\bar  \nabla_\Sigma \tilde f (\pi_{\pa E}(x))|. 
\end{split}
\]
We use  Lemma \ref{lem:Leibniz} with $k =2$ for $T_1 = T_2 = B_\Sigma$  and $T_3 = \bar  \nabla_\Sigma \tilde f $, and   \eqref{eq:geometric-regularity-33}  to infer
\[
\|\nabla^{3} (\tilde f\circ \pi_{\pa E})\|_{L^2(\Gamma)} \leq C( \| \tilde f \|_{H^3(\Sigma)} + C(1+ \|B_\Sigma\|_{H^2(\Sigma)}) \| \bar \nabla \tilde f \|_{L^\infty(\Sigma)}).
\]
We use  Sobolev embedding, interpolation  in Proposition \ref{prop:interpolation} and  \eqref{eq:geometric-regularity-3} 
\[
 \| \bar \nabla \tilde f \|_{L^\infty(\Sigma)} \leq C  \|  \tilde f \|_{W^{2,4}(\Sigma)} \leq C  \| \tilde f \|_{H^3(\Sigma)}^{\frac12} \| \tilde f \|_{H^2(\Sigma)}^{\frac12} \leq Ch^{-\frac18}. 
\]
Therefore  the two inequalities  above together with    \eqref{eq:geometric-regularity-3} and \eqref{eq:upgr-reg-0}  yield
\begin{equation}\label{eq:geometric-regularity-8}
\|\nabla^{3} (\tilde f\circ \pi_{\pa E})\|_{L^2(\Gamma)} \leq C(h^{-\frac14} +   K_0 h^{-\frac18}) \leq Ch^{-\frac14}
\end{equation}
for $h \leq h_0 = h_0(r_0,R_0,K_0,\delta)$.

Recall the formula \eqref{eq:geometric-laplacian} for  the Laplacian of $ \tilde f\circ \pi_{\pa E}$  on $\pa F$, and  that by  the first equation in  \eqref{eq:Euler-Lagrange-4} it holds $\nabla_{\pa F} H_F = - \nabla_{\pa F} (\tilde f\circ \pi_{\pa E}) $.  Therefore we may estimate point-wise on $\pa F $ using \eqref{eq:geometric-regularity-33}
\[
|\nabla_{\pa F}\Delta_{\pa F} (\tilde f\circ \pi_{\pa E}) | \leq C\big(|\nabla^{3} (\tilde f\circ \pi_{\pa E})|  + |\nabla^2 (\tilde f\circ \pi_{\pa E})| + |\nabla (\tilde f\circ \pi_{\pa E})|  + |\nabla (\tilde f\circ \pi_{\pa E})|^2\big). 
\]
We infer  by the Sobolev embedding,  by \eqref{eq:geometric-regularity-7} and  \eqref{eq:geometric-regularity-8}  that 
\[
\begin{split}
\|\nabla_{\pa F}\Delta_{\pa F} (\tilde f\circ \pi_{\pa E}) \|_{L^2(\pa F)} &\leq C \|\nabla^{3} (\tilde f\circ \pi_{\pa E})\|_{L^2(\pa F)}  +C\| \nabla^{2} (\tilde f\circ \pi_{\pa E}) \|_{L^2(\pa F)}   \\
 &\,\,\,\,\,\,\,\,\,\,+C\| \nabla (\tilde f\circ \pi_{\pa E}) \|_{L^2(\pa F)}  + C\|\nabla (\tilde f\circ \pi_{\pa E}) \|_{L^4(\pa F)}^2 \\ 
&\leq  C \|\nabla^{3} (\tilde f\circ \pi_{\pa E})\|_{L^2(\pa F)}   +C \| \tilde f\circ \pi_{\pa E} \|_{H^2(\Gamma)}  + C  \| \tilde f\circ \pi_{\pa E} \|_{H^2(\Gamma)}^2\\
&\leq C h^{-\frac14} +C (1+ K_0)\leq C h^{-\frac14},
\end{split}
\]
for $h \leq h_0 = h_0(r_0,R_0,K_0)$.   This concludes the proof.

\end{proof}

\section{Proof of the main theorems}

\subsection{Iteration Lemma}

%

We  consider three sets $E, F $ and $G$  such that  $E \subset \R^3$ satisfies the assumptions \eqref{eq:assumptions-E} and  $\|\Delta_{\pa E} H_E\|_{H^1(\pa E)} \leq K_0 h^{-\frac14}$, and  $F$ is a minimizer of \eqref{def:min-prob}. It follows from
 Lemma \ref{lem:upgrade-regularity} and Lemma \ref{lem:geometric-regularity}   that also $F$  satisfies the assumptions \eqref{eq:assumptions-E} with different  constants, i.e.,   it satisfies the UBC with radius $r_0/2$, 
\[
\|\Delta_{\pa F} H_F\|_{L^2(\pa F)} \leq K_1  \quad \text{and} \quad   \|\Delta_{\pa F} H_F\|_{H^1(\pa F)} \leq K_1  h^{-\frac14}
\]
for some $K_1$ depending on $r_0,R_0,K_0$  and $F\subset B_{2R_0}$.  We then proceed and consider a minimizer $G$ of the problem 
\begin{equation}\label{eq:iteration-1}
\min   \Big{\{} P(G)+\frac{d_{H^{-1}}(G;F)^2}{2h} :   G \Delta F \subset \overline{\mathcal{N}_{\delta}}(\pa F) \Big{\}} . 
\end{equation}
We may again apply Proposition \ref{prop:apirori-estimates},  Lemma \ref{lem:upgrade-regularity}  and Lemma \ref{lem:geometric-regularity} to deduce that when $\delta \leq \delta_3$, where $\delta_3$ depends on $r_0,K_1, R_0$  and $\delta_2$ from Proposition \ref{prop:apirori-estimates} and $h \leq  h_1$,   we may write the boundary of $G$ as a normal deformation with a heightfunction  $\psi_F :\pa F \to \R$, 
\[
\pa G = \{ x + \psi_{F}(x)\nu_F(x) : x \in \pa F\}
\]
and $\|\psi_F\|_{L^2(\pa F)}\leq C_3h$, $\|\psi_F\|_{H^4(\Gamma)}\leq C_3$, where $\Gamma = \pa F$. In order to avoid confusion we denote the heighfunction of $\pa F$ over $\pa E$ as $\psi_E$, i.e., 
\[
\pa F = \{ x + \psi_{E}(x)\nu_E(x) : x \in \pa E\}. 
\]
Since $K_1$ depends on $r_0, R_0$ and $K_0$, then eventually $\delta_3$ also depends on $r_0,R_0, K_0$ and $\delta_2$. 

\begin{lemma}\label{lem:iteration}
Assume $E \subset \R^3$ is $C^5$-regular, satisfies the assumptions \eqref{eq:assumptions-E} and  $\|\Delta_{\pa E} H_E\|_{H^1(\pa E)} \leq K_0 h^{-\frac14}$. Assume  $F\subset \R^3$ is a minimizer of \eqref{def:min-prob} with $\delta < \delta_3$, where $\delta_3$ depends on $r_0,R_0, K_0$  and $\delta_2 $ from Proposition \ref{prop:apirori-estimates}, and assume $G\subset \R^3$  is a minimizer of \eqref{eq:iteration-1} and denote the associated  heightfunctions by $\psi_E, \psi_F$ and let $\xi_{F,E}, \xi_{G,F}$ as in \eqref{def:normal-velo}.  Then, there exists a constant $C_4$, depending on $r_0, R_0$ and $K_0$, such that
\begin{equation}\label{eq:iteration-estimate}
\int_{\pa F}\Big(\xi_{G,F}^2+\frac{h}{2}|\Delta_{\pa F}\psi_{F}|^2\Big)\,d\H^2\leq(1+C_4h)\int_{\pa E}\Big( \xi_{F,E}^2 +\frac{h}{8}|\Delta_{\pa E}\psi_{E}|^2 \Big)\,d\H^2,
\end{equation}
for $h \leq h_1$, where $h_1 = h_1(r_0,R_0,K_0, \delta)$. 
\end{lemma}
\begin{proof}
By the discussion at the beginning of the section, we may deduce that  Proposition \ref{prop:apirori-estimates},  Lemma \ref{lem:upgrade-regularity}  and Lemma \ref{lem:geometric-regularity} hold also for $F$ and $G$ with possiby smaller uniform ball radius  and smaller treshhold for $\delta$ which we denote by $\delta_3 \leq \delta_2$.  In particular, we may write the Euler-Lagrange equation \eqref{eq:Euler-Lagrange-2} for $G$, which reads as  
\[
\frac{1}{h}\xi_{G,F} = \Delta_{\pa F} \left( - \Delta_{\pa F}\psi_F  + H_F  \right) +  \Delta_{\pa F} R_0 \qquad \text{on }\, \pa F,
\]
where the error term is of the form  \eqref{eq:error-R-0}. We multiply this by $\xi_{G,F} $ and integrate by parts to infer
\begin{equation}\label{eq:itetation-2}
\frac{1}{h}\int_{\pa F} \xi_{G,F}^2\,d\H^2 =   \int_{\pa F}  \left( - \Delta_{\pa F}\psi_F  + H_F  \right) (\Delta_{\pa F} \xi_{G,F})\,d\H^2 + \int_{\pa F} R_0 (\Delta_{\pa F} \xi_{G,F})\,d\H^2.
\end{equation}
Recall that by  \eqref{def:normal-velo}
\[
\xi_{G,F} = \psi_F + \frac{H_F}{2} \psi_F^2 + \frac{K_F}{3} \psi_F^3. 
\]
Arguing as in \eqref{eq:upgr-reg-88} we deduce that for every $\eps>0$ it holds  by the smallness of $ \|\psi_F\|_{L^\infty(\pa F)}$ that
\[
\| \Delta_{\pa F}  (\frac{H_F}{2} \psi_F^2)  \|_{L^2(\pa F)} \leq C  \|\psi_F\|_{L^\infty(\pa F)} \|\psi_F\|_{H^2(\Gamma) }\leq \eps  \|\psi_F\|_{H^2(\Gamma) }
\]
for  $\Gamma = \pa F$ when $h$ is small. A similar argument yields
\[
\| \Delta_{\pa F}  (\frac{K_F}{3} \psi_F^3)  \|_{L^2(\pa F)} \leq C  \|\psi_F\|_{L^\infty(\pa F)} \|\psi_F\|_{H^2(\Gamma) } \leq \eps  \|\psi_F\|_{H^2(\Gamma) }.
\]
Therefore we have by \eqref{eq:itetation-2} using the above and  \eqref{eq:hilbert-norm} 
\begin{equation}\label{eq:iteration-3}
\begin{split}
\frac{1}{h}&\int_{\pa F} \xi_{G,F}^2\,d\H^2  +\int_{\pa F} ( \Delta_{\pa F}\psi_F)^2\,d\H^2 \\
&\leq    \int_{\pa F}  H_F (\Delta_{\pa F} \xi_{G,F})\,d\H^2 + \eps  \|\Delta_{\pa F}\psi_F\|_{L^2(\pa F)}^2  +C\eps \|\psi_F\|_{L^2(\pa F)}^2\\
&\,\,\,\,\,\,\,\,\,\,+    C ( \|\Delta_{\pa F}\psi_F\|_{L^2(\pa F)}  + \|\psi_F\|_{L^2(\pa F)})\|R_0\|_{L^2(\pa F)}. 
\end{split}
\end{equation}

We need to bound $\|R_0\|_{L^2(\pa F)}$. We use the form in  \eqref{eq:error-R-0} to bound the error term point-wise on $\Gamma = \pa F$ as
\[
|R_0| \leq C(|\psi_F| + |\bar \nabla \psi_F|) (1+ |\bar \nabla^2 \psi_F| + |\bar \nabla (\psi_F \, B_\Gamma)|). 
\]
By  Lemma \ref{lem:geometric-regularity}, Lemma \ref{lem:hilbert-norm} and by Sobolev embedding it holds 
\begin{equation}\label{eq:iteration-33}
\|B_\Gamma\|_{L^\infty(\Gamma)}, \|B_\Gamma\|_{H^2(\Gamma)}  \leq C \quad \text{and} \quad \|\bar \nabla B_\Gamma\|_{L^\infty(\Gamma)} \leq Ch^{-\frac14}. 
\end{equation}
Moreover Proposition \ref{prop:apirori-estimates}  with $\eps = \frac14$ implies $\|\psi_F\|_{L^\infty(\pa F)} \leq C h^{\frac14}$.  Therefore we have 
\[
 |\bar \nabla (\psi_F \, B_\Gamma)| \leq  |\bar \nabla \psi_F| |B_\Gamma|  +|\psi_F || \bar \nabla B_\Gamma| \leq C. 
\]
Thus it holds by interpolation from Proposition \ref{prop:interpolation}, by the smallness of $\|\psi_F\|_{C^1}$ and by  \eqref{eq:hilbert-norm}
\[
\begin{split}
\|R_0\|_{L^2(\pa F)} &\leq C\|\psi_F\|_{C^1} \|\psi_F\|_{H^2(\Gamma) } + C\|\psi_F\|_{H^1(\Gamma) }\leq \eps  \|\psi_F\|_{H^2(\Gamma) } + C_\eps \|\psi_F\|_{L^2(\pa F) }\\
&\leq  \eps  \|\Delta_{\pa F }\psi_F\|_{L^2(\pa F) } + C_\eps \|\psi_F\|_{L^2(\pa F) }.
\end{split}
\]
Recall also that \eqref{eq:geometric-44} implies $|\psi_F(x)| \leq 2| \xi_{G,F}(x)|$ for all $x \in \pa F$.  Hence we may estimate  \eqref{eq:iteration-3}  as
\begin{equation}\label{eq:iteration-5}
\left(\frac{1}{h} - C\right)\int_{\pa F} \xi_{G,F}^2\,d\H^2  + \frac12\int_{\pa F} ( \Delta_{\pa F}\psi_F)^2\,d\H^2 \leq   \int_{\pa F}  H_F (\Delta_{\pa F} \xi_{G,F})\,d\H^2 
\end{equation}
when $\eps$ is chosen small enough.

Let us next analyze  the right-hand-side of \eqref{eq:iteration-5}. We integrate by parts and use a change of variables with the diffeomorphism $\Psi_E : \pa E \to \pa F$, $\Psi_E(x) = x + \psi_E(x) \nu_E(x)$
\[
\begin{split}
\int_{\pa F}  H_F (\Delta_{\pa F} \xi_{G,F})\,d\H^2  &=  - \int_{\pa F} ( \nabla_{\pa F}  H_F \cdot  \nabla_{\pa F}   \xi_{G,F}  ) \,d\H^2 \\
&=-\int_{\pa E}\big( ( \nabla_{\pa F}  H_F \cdot  \nabla_{\pa F}   \xi_{G,F}  ) \circ \Psi_E  \big) \, J_{\pa E}\Psi_E  \,d\H^2.
\end{split}
\]
We use  \eqref{eq:change-variables-function-2} for the function $f = H_F \circ \Psi_E : \pa E \to \R$. Notice that for this function it holds $f \circ \pi_{\pa E} = H_F$ on $\pa F$. We deduce by   \eqref{eq:change-variables-function-2} that  for $x \in \pa E$ it holds 
\[
(\nabla_{\pa F}   H_F)(\Psi_E(x))  = (I + A(\psi_E(x), \nabla_{\pa E} \psi_E(x), \nu_E(x),B_E(x)))\nabla_{\pa E}( H_F \circ \Psi_E)(x),
\]
for a matrix which satisfies $ A(0,0, \cdot, \cdot ) = 0$. Similarly it holds 
\[
(\nabla_{\pa F}   \xi_{G,F})(\Psi_E(x))  = (I + A(\psi_E(x),\nabla_{\pa E} \psi_E(x), \nu_E(x),B_E(x)))\nabla_{\pa E}(  \xi_{G,F} \circ \Psi_E)(x). 
\]
Therefore by the above and by the formula for the Jacobian \eqref{eq:jacobian-psi} we have
\begin{equation}\label{eq:iteration-6}
\int_{\pa F}  H_F (\Delta_{\pa F} \xi_{G,F})\,d\H^2  \leq  -\int_{\pa E}\big( ( \nabla_{\pa E}  (H_F \circ \Psi_E) \cdot  \nabla_{\pa E}  ( \xi_{G,F} \circ \Psi_E  )  \big) \, \sqrt{J_{\pa E}\Psi_E}  \,d\H^2 + \|R_1\|_{L^1(\pa E)},
\end{equation}
where the error term satisfies the point-wise bound $x \in \pa E$ 
\[
|R_1(x)| \leq C(|\psi_E(x)| + |\nabla_{\pa E} \psi_E(x)|) |\nabla_{\pa F} H_F(\Psi_E(x))| \,  |\nabla_{\pa F} \xi_{G,F}(\Psi_E(x))| .
\]
By H\"older's inequality 
\[
\|R_1\|_{L^1(\pa E)} \leq C\|\psi_E\|_{W^{1,4}(\pa E)} \| \nabla_{\pa F} H_F \|_{L^{4}(\pa F)} \|  \nabla_{\pa F} \xi_{G,F}\|_{L^2(\pa F)}. 
\]
The estimate \eqref{eq:geometric-44} implies $ \|  \nabla_{\pa F} \xi_{G,F}\|_{L^2(\pa F)} \leq C \|\psi_F\|_{H^1(\pa F)}$, while \eqref{eq:iteration-33}  and Sobolev embedding give
$\| \nabla_{\pa F} H_F \|_{L^{4}(\pa F)} \leq \| B_\Gamma  \|_{H^{2}(\Gamma)} \leq C$. Therefore we may bound by interpolation from Proposition \ref{prop:interpolation}, by  \eqref{eq:hilbert-norm} and by  \eqref{eq:geometric-44}
\begin{equation}\label{eq:iteration-66}
\begin{split}
\|R_1\|_{L^1(\pa E)} &\leq C\|\psi_E\|_{W^{1,4}(\pa E)}\|\psi_F\|_{H^1(\pa F)}\\
&\leq \eps \|\Delta_{\pa E} \psi_E\|_{L^2(\pa E)}^2  + C_\eps  \|  \xi_{F,E}\|_{L^2(\pa E)}^2+  \eps \|\Delta_{\pa F} \psi_F\|_{L^2(\pa F)}^2  + C_\eps  \|  \xi_{G,F}\|_{L^2(\pa F)}^2.
\end{split}
\end{equation}

We integrate \eqref{eq:iteration-6} by parts and have 
\begin{equation}\label{eq:iteration-7}
\begin{split}
\int_{\pa F}  H_F (\Delta_{\pa F} \xi_{G,F})\,d\H^2  \leq   \int_{\pa E} &\big( \Delta_{\pa E}  (H_F \circ \Psi_E ) \big)  ( \xi_{G,F} \circ \Psi_E  )   \, \sqrt{J_{\pa E}\Psi_E}  \,d\H^2 \\
&+ \|R_1\|_{L^1(\pa E)}+ \|R_2\|_{L^1(\pa E)},
\end{split}
\end{equation}
where the first error term satisfies \eqref{eq:iteration-66} and  the second error term satisfies the following  point-wise bound on $\Sigma = \pa E$ by the formula for $J_{\pa E}\Psi_E$ in   \eqref{eq:jacobian-psi}
\[
|R_2(x)| \leq C |\psi_F(\Psi_E(x))| |\nabla_{\pa F} H_F(\Psi_E(x))|  (|\psi_E(x)| + |\nabla_{\pa E} \psi_E(x)|) (1+  |\bar \nabla_{\Sigma}^2 \psi_E(x)| + |\bar \nabla_{\Sigma} B_\Sigma(x)|) .
\]
By H\"older's inequality 
\[
\|R_2\|_{L^1(\pa E)} \leq C \|\psi_F\|_{L^\infty(\pa F)} \| \nabla_{\pa F} H_F \|_{L^{4}(\pa F)} \|\psi_E\|_{W^{1,4}(\pa E)} (1 + \|  \psi_E\|_{H^2(\Sigma)} + \|  B_\Sigma\|_{H^1(\Sigma)}).  
\]
By  \eqref{eq:iteration-33}  it holds  $ \| \nabla_{\pa F} H_F \|_{L^{4}(\pa F)}\leq \|B_\Gamma\|_{H^2(\Gamma)} \leq C$ and Lemma \ref{lem:upgrade-regularity}  implies $\|  \psi_E\|_{H^2(\Sigma)} \leq C$. We use  the inequality  \eqref{eq:hilbert-norm},  Sobolev embedding and interpolation from Proposition \ref{prop:interpolation} to get 
\begin{equation}\label{eq:iteration-8}
\begin{split}
\|R_2\|_{L^1(\pa E)} & \leq C\big(\|\psi_E\|_{W^{1,4}(\pa E)} ^2+ \|\psi_F\|_{L^\infty(\pa F)}^2\big) \\
&\leq   \eps \|\Delta_{\pa E} \psi_E\|_{L^2(\pa E)}^2  + C_\eps  \|  \xi_{F,E}\|_{L^2(\pa E)}^2+  \eps \|\Delta_{\pa F} \psi_F\|_{L^2(\pa F)}^2  + C_\eps  \|  \xi_{G,F}\|_{L^2(\pa F)}^2.
\end{split}
\end{equation}


Since $F$ is a minimizer it satisfies the Euler-Lagrange equation \eqref{eq:Euler-Lagrange-3}, which we write as a one equation on $\pa E$
\[
\frac{1}{h} \xi_{F,E}  = \Delta_{\pa E} ( H_F \circ \Psi_E) \qquad \text{on } \, \pa E.   
\] 
We use this in \eqref{eq:iteration-7} and   have by Young's inequality and change of variables  
\[
\begin{split}
\int_{\pa F}  &H_F (\Delta_{\pa F} \xi_{G,F})\,d\H^2 \\
 &\leq  \frac{1}{h} \int_{\pa E} \xi_{F,E}  \,  ( \xi_{G,F} \circ \Psi_E  )   \, \sqrt{J_{\pa E}\Psi_E}  \,d\H^2  + \|R_1\|_{L^1(\pa E)} + \|R_2\|_{L^1(\pa E)} \\
&\leq \frac{1}{2h} \int_{\pa E} ( \xi_{G,F} \circ \Psi_E  )^2 \, J_{\pa E}\Psi_E   \,d\H^2 + \frac{1}{2h} \int_{\pa E}   \xi_{F,E}^2 \,d\H^2 + \|R_1\|_{L^1(\pa E)} + \|R_2\|_{L^1(\pa E)} \\
&=  \frac{1}{2h} \int_{\pa F}  \xi_{G,F}^2   \,d\H^2 + \frac{1}{2h} \int_{\pa E}   \xi_{F,E}^2 \,d\H^2 + \|R_1\|_{L^1(\pa E)} + \|R_2\|_{L^1(\pa E)} .
\end{split}
\]
The above,  \eqref{eq:iteration-5}, \eqref{eq:iteration-66} and \eqref{eq:iteration-8}    yield 
\[
\begin{split}
\left(\frac{1}{2h} - C\right) \int_{\pa F} \xi_{G,F}^2 \,d\H^2&+  \frac13\int_{\pa F} |\Delta_{\pa F}\psi_{F}|^2 \,d\H^2\\
&\leq  \left(\frac{1}{2h} + C\right) \int_{\pa E}\xi_{F,E}^2 \,d\H^2    + \frac{1}{20} \int_{\pa E} |\Delta_{\pa E}\psi_{E}|^2 \,d\H^2.
\end{split}
\]
when $\eps$ is chosen small enough, and the claim follows. 

\end{proof}

\subsection{Proof of Theorem \ref{thm:thm2} and Theorem \ref{thm:thm1} }

We remark that under the assumptions  $\|\psi\|_{L^2}\leq L_0h$ and $\|\psi\|_{H^4}\leq L_0$ in Lemma \ref{lem:geometric-regularity},  when $K_0$ in the assumption   \eqref{eq:assumptions-E}   is chosen large enough,  we know that the minimizer $F$ satisfies 
\[
\|\Delta_{\pa F} H_{F} \|_{L^2(\pa F)}  \leq K_0 \qquad \text{and} \qquad  \|\Delta_{\pa F} H_{F} \|_{H^1(\pa F)} \leq K_0 h^{-\frac14}.
\]
This means that we are able to keep the a priori regularity estimates assumed in \eqref{eq:assumptions-E} except possibly the radius for the uniform ball condition which may decrease from $r_0$ to $r_0/2$. 
We may ensure that we also keep the uniform ball condition  using the following lemma. 
 
\begin{lemma}
\label{lem:keep-UBC}
Assume that $E_0 \subset \R^3$ satisfies the UBC with radius $r_0$, $E_0 \subset B_{R_0}$ and is uniformly $C^{3}$-regular. Assume further that $F\subset \R^n$ is a $(\Lambda,\alpha)$-minimizer of the perimeter according to Definition \ref{def:lambda-mini} and 
\[
\|B_F\|_{L^\infty(\pa F)} \leq K \quad \text{and} \quad \|\Delta_{\pa F} H_F\|_{L^2(\pa F)} \leq K,
\]
for some $K>0$. Then there is $\delta$, depending on $r_0,R_0, \Lambda$, $\alpha$ and $K$, such that if $F \Delta E_0 \subset \mathcal{N}_\delta(\pa E_0)$ then  $F$ satisfies the UBC with radius $r_0/2$. 
 \end{lemma} 
\begin{proof}
Lemma \ref{lem:omega-mini} implies that when $\delta$ is small then there is $\psi \in C^{1,\gamma}(\pa E_0)$ such that 
\[
\pa F = \{ x + \psi(x)\nu_{E_0}(x) : x \in \pa E_0\}
\]
with $\|\psi\|_{C^{1,\gamma}(\pa E_0)}$ small. Note that  from  $F \Delta E_0 \subset \mathcal{N}_\delta(\pa E_0)$ it follows  $\|\psi\|_{L^\infty} \leq \delta$.  The assumptions on $F$ together with Lemma \ref{lem:hilbert-norm} and Sobolev embedding yield for all $\beta\in[0,1)$, $\|B_{\Gamma}\|_{C^{0,\beta}(\Gamma)} \leq C$, where $\Gamma=\pa F$. By straightforward calculation we then deduce that $\|\psi\|_{C^{2,\beta}(\Gamma)} \leq C$. By the interpolation inequality \eqref{eq:Holder-interpolate} we then have $\|\psi\|_{C^{2}(\Gamma)} \leq C\delta^{\frac{\beta}{2+\beta}}$. We deduce by \cite[Theorem~2.6]{Dal} and \cite[Remark 3.4.6]{Antonia} that $F$ satisfies the UBC with radius $r_0/2$ when $\delta$ is small. 
\end{proof}

\begin{proof}[\textbf{Proof of Theorem \ref{thm:thm2}}]
Assume $E_0 \subset B_{R_0/2} \subset \R^3$ is a $C^5$-regular set which satisfies the uniform ball condition with radius $r = 2r_0 $ and $\|\Delta_{\pa E_0} H_{E_0} \|_{L^2(\pa E_0)} \leq K$. We fix a large constant $K_0 = K_0(r_0,R_0,K)$, whose choice will be clear later. Let us then choose $\delta'$, which is from Lemma \ref{lem:keep-UBC}, and    $\delta < \min\{\delta_3, \delta'\}$, where $\delta_3$ is from Lemma  \ref{lem:iteration}. We consider an approximate sequence $E_k^{h,\delta}$  minimizing the energy \eqref{def:min-prob} and simplify the notation by  setting $E_k = E_k^{h,\delta}$. We obtain from Proposition \ref{prop:apirori-estimates}, Lemma \ref{lem:upgrade-regularity} and Lemma \ref{lem:geometric-regularity} that we may write 
\begin{equation}\label{eq:thm-2-height}
\pa E_1  = \{ x + \psi_0(x) \nu_{E_0}(x) : x \in \pa E_0 \}
\end{equation}
and it holds  $\|\psi_0 \|_{L^2(\pa E_0)} \leq L_0 h$ and $\|\psi_0 \|_{H^4(\Sigma_0)} \leq L_0$, for $\Sigma_0 = \pa E_0$ and  $L_0= L_0(r_0, R_0,K)$. 

We thus have obtained that $E_0$ and $\psi_0$ satisfy by the above and by interpolation from Proposition \ref{prop:interpolation}
\[
\|\psi_0 \|_{L^2(\pa E_0)} \leq L_0 h \quad \text{and} \quad \|\Delta_{\pa E_0}\psi_0 \|_{L^2(\Sigma_0)} \leq  L_0\sqrt{h},
\]
 by possibly increasing the original constant $L_0$, and 
\[
\|\Delta_{\pa E_0} H_{E_0} \|_{L^2(\pa E_0)} \leq  K_0 \quad \text{and} \quad \|\Delta_{\pa E_0} H_{E_0} \|_{H^1(\pa E_0)} \leq  K_0h^{-\frac14}. 
\]
 We  deduce by Lemma  \ref{lem:geometric-regularity} that there is a constant $C_2 = C_2(r_0,R_0, L_0)$ such that  
\[
\|\Delta_{\pa E_1} H_{E_1} \|_{L^2(\pa E_1)} \leq C_2(1+ \sqrt{K_0}) \leq K_0 \quad \text{and} \quad \|\Delta_{\pa E_1} H_{E_1} \|_{H^1(\pa E_1)} \leq C_2h^{-\frac14} \leq K_0h^{-\frac14},
\]
when $K_0$ is chosen  large enough, and that $E_1$ satisfies the UBC with radius $r_0$. 

We denote by $k_0 \in \mathbb{N}$  the largest   index  such that $\sigma_0 = k_0 h \leq 1$ and   for all $k \leq k_0$ it holds 
\[
\pa E_k \subset \mathcal{N}_{\delta}(\pa E_0).
\]
We claim that for every $k \leq k_0$, $E_k$ satisfies the UBC with radius $r_0$ and  it holds 
\begin{equation}\label{eq:induction-reg}
  \|\Delta_{\pa E_k} H_{E_k}\|_{L^2(\pa E_k)} \leq K_0 \quad \text{and} \quad \|\Delta_{\pa E_k} H_{E_k}\|_{H^1(\pa E_k)} \leq K_0 h^{-\frac14}.
\end{equation}

We prove \eqref{eq:induction-reg} by induction.  Recall that it holds for $k=1$, and assume it holds for $k-1$.  By applying again Proposition \ref{prop:apirori-estimates}, Lemma \ref{lem:upgrade-regularity} and Lemma \ref{lem:geometric-regularity}  we may write the set $E_k$ as 
\[
\pa E_k  = \{ x + \psi_{k-1}(x) \nu_{E_{k-1}}(x) : x \in \pa E_{k-1} \}
\]
with $\|\psi_{k-1} \|_{L^2(\pa E_{k-1})} \leq C_1 h$ and $\|\psi_{k-1} \|_{H^4(\Sigma_{k-1})} \leq C_1$ for $C_1 = C_1(r_0,R_0, K_0)$, where $\Sigma_{k-1}=\pa E_{k-1}$. We denote $\xi_{k-1} = \xi_{E_{k}, E_{k-1}}$, the function defined  in \eqref{def:normal-velo}, and  have by Lemma \ref{lem:iteration}
\[
\int_{\pa E_{j}}\Big(\xi_{j}^2+\frac{h}{2}|\Delta_{\pa E_j}\psi_{j}|^2\Big)\,d\H^2\leq(1+C_4h)\int_{\pa E_{j-1}}\Big( \xi_{j-1}^2 +\frac{h}{8}|\Delta_{\pa E_{j-1}}\psi_{j-1}|^2 \Big)\,d\H^2
\]
for every $j \leq k$.  By iterating the above estimate and using  $\|\xi_0 \|_{L^2(\pa E_0)} \leq 2\|\psi_0 \|_{L^2(\pa E_0)} \leq 2 L_0 h$ and  $\|\Delta_{\pa E_0} \psi_0 \|_{L^2(\pa E_0)} \leq L_0 \sqrt{h}$,  we obtain  for a constant  $C_4 = C_4(r_0,R_0, K_0)$
\[
\begin{split}
\int_{\pa E_{k}}\Big(\xi_{k}^2 +\frac{h}{4} \sum_{j=1}^k |\Delta_{\pa E_j}\psi_{j}|^2\Big)\,d\H^2 &\leq (1+C_4h)^k \int_{\pa E_{0}}\Big( \xi_{0}^2 +\frac{h}{8}|\Delta_{\pa E_{0}}\psi_{0}|^2 \Big)\,d\H^2\\
&\leq 5e^{C_4kh} L_0^2 h^2\\
&\leq5e^{C_4k_0h} L_0^2 h^2  \leq 10 L_0^2 h^2, 
\end{split}
\]
provided that $k_0 h = \sigma_0 = \sigma_0(r_0, R_0, K_0)$ is small. Therefore, by possibly  increasing the value of $L_0$ we deduce by the above and   \eqref{eq:geometric-44}
\begin{equation}\label{eq:thm-2-1}
\| \psi_{k}\|_{L^2(\pa E_k)}^2 + h  \sum_{j=0}^k  \| \Delta_{\pa E_j} \psi_{j}\|_{L^2(\pa E_j)}^2   \leq L_0^2 h^2.
\end{equation}
Lemma \ref {lem:geometric-regularity} implies that it holds for $C_2 = C_2(r_0,R_0,L_0)$
\[
\|\Delta_{\pa E_k}H_{E_k}\|_{L^2(\pa E_k)} \leq  C_2(1 + \sqrt{K_0}) \leq K_0  \qquad \text{and} \qquad  \|\Delta_{\pa E_k}H_{E_k}\|_{H^1(\pa E_k)} \leq C_2 h^{-\frac14}\leq K_0  h^{-\frac14},
\]
when  $K_0$ is large. Finally we deduce by Lemma \ref{lem:keep-UBC} that the set $E_k$ satisfies the UBC with radius $r_0$, when $\delta'$ is chosen small enough. This proves \eqref{eq:induction-reg}. 

Let us then show that if we denote $\sigma_0 = k_0h$, then $\sigma_0$ is uniformly bounded from below. To this aim we assume that for the set $E_{k_0}$ there is  a point $x \in \pa E_{k_0}$ such that $|d_{E_0}(x)| \geq \delta/2$. 
Proposition \ref{prop:apirori-estimates} implies that $E_{k_0}$ is a $(\Lambda,\alpha)$-minimizer of the perimeter for any $\alpha \in (0,\frac13)$, where $\Lambda$ depends on $r_0, R_0, K_0$ and $ \alpha$. Then by  $|d_{E_0}(x)| \geq \delta/2$ and by  standard density estimate 
we  have 
\begin{equation}\label{eq:thm-2-2}
|E_{k_0} \Delta E_0| \geq c \delta^3. 
\end{equation}

Let $k \leq k_0$. Following the calculations in the proof of Lemma \ref{lem:distance-finite} and using \eqref{eq:geometric-44}  we may write 
\[
|E_{k} \Delta E_{k-1}| = \int_{\pa E_{k-1}} |\xi_{k-1}|\, d \mathcal{H}^2  \leq 2  \int_{\pa E_{k-1}} |\psi_{k-1}|\, d \mathcal{H}^2 . 
\]
Then we may estimate by  \eqref{eq:thm-2-1} 
\[
\begin{split}
|E_{k_0} \Delta E_0|  &\leq \sum_{j=1}^{k_0}|E_{k} \Delta E_{k-1}|  \leq 2   \sum_{j=0}^{k_0-1} \|\psi_{j}\|_{L^1(\pa E_{j})} \\
&\leq C  \sum_{j=0}^{k_0-1} \|\psi_{j}\|_{L^2(\pa E_{j})} \leq C k_0 h =  C \sigma_0.  
\end{split}
\]
The inequality \eqref{eq:thm-2-2}  proves that $\sigma_0$ is uniformly bounded from below. 

We have thus proved a uniform $H^4$-regularity estimate \eqref{eq:induction-reg} for  the sets $E_k$ for $k \leq k_0$ where $k_0 h = \sigma_0$ is uniformly bounded from below.   Again Proposition \ref{prop:apirori-estimates}  implies that   the sets $E_k$ are $(\Lambda, \alpha)$-minimizers of the perimeter for constants that are independent of $h$. Therefore by Lemma \ref{lem:omega-mini} we deduce that we may write  
\[
\pa E_k = \{ x + u_k(x)\nu_{E_0}(x) : x \in \pa E_0 \} \quad \text{for }\, u_k : \pa E_0 =  \Sigma_0: \to \R,
\]
when $\sigma_0$ is sufficiently small, $\|u_k\|_{C^{1,\gamma}(\Sigma_0)} \leq \delta^\gamma$ for some $ \gamma>0$,  and by the uniform bounds from \eqref{eq:induction-reg} we deduce that $\|u_k\|_{H^4(\Sigma_0)} \leq C$. In particular, $u_k$ are uniformly $C^{2,\alpha}$-regular.  Let us define the discrete normal velocity on $\pa E_k$ as  $v_k : \pa E_k \to \R$,  
\[
v_k(y) := \frac{\psi_{k}(y)}{h}.
\]
Define $\Phi_k: \pa E_0 \to \pa E_k$, $\Phi_k(x) = x + u_k(x)\nu_{E_0}(x)$ and denote its tangential Jacobian as $J_{\pa E_0}\Phi_k$.  We claim that it holds 
\begin{equation}
\label{eq:normal-velo-disc}
\Big{\|}v_k\circ \Phi_k - \frac{u_{k+1} - u_k}{|N_k| \, h } \Big{\|}_{L^2(\pa E_0)}\leq \eps(h)
\end{equation}
for $\eps(h) \to 0$ as $h \to 0$, where
\[
|N_k(x)| = \frac{J_{\pa E_0}\Phi_k(x)}{1 + H_{E_0}u_k(x) + K_{E_0}u_k^2(x)}. 
\]
The notation is justified by the fact that $|N_k(x)|$ defined above is in fact the norm of the vector  $N_{k}(x)$ defined in \eqref{eq:meancurva-exp-0} for $E = E_0$ and $F = E_k$, which is not obvious but follows from standard calculation.  
We note that by the estimates from the previous section   we have $\|\psi_k \circ \Phi_k \|_{C^1} \leq  h^{\gamma}$  for $\gamma >0$ and $\|\psi_k \circ \Phi_k \|_{L^2} \leq C h$. Moreover, as  $\|u_k\|_{C^{1,\gamma}(\Sigma_0)} \leq \delta^\gamma$ then it is easy to see that 
$|u_{k+1}(x)- u_{k}(x)| \leq C|\psi_k(\Phi_k(x))|$ for all $x \in \pa E_0$ and $\|u_{k+1}- u_{k}\|_{C^1}\leq h^\gamma$ by possibly decreasing the value of $\gamma$ if needed.  We will use these repeatedly below. 

In order to show \eqref{eq:normal-velo-disc} we fix  $f: \pa E_0 \to \R$ such that $\|f\|_{C^1}\leq C h^\gamma$, define $\Phi_\tau: \pa E_0 \to \R^3$ as $\Phi_\tau(x) = x + \tau \nu_{E_0}(x)$, and calculate the following integral as in the proof of  Lemma \ref{lem:distance-finite}
\[
\begin{split}
\int_{\R^3} &f(\pi_{\pa E_0}(x))(\chi_{E_{k+1}} - \chi_{E_{k}})\, dx \\
&= \int_{\pa E_0} f(x) \int_{-r_0}^{r_0} (\chi_{E_{k+1}}(\Phi_\tau(x)) - \chi_{E_{k}}(\Phi_\tau(x)) J_{\pa E_0}\Phi_\tau(x) \, d \tau d \H^2(x)\\
&= \int_{\pa E_0} f(x) \int_{u_k(x)}^{u_{k+1}(x)} J_{\pa E_0}\Phi_\tau(x) \, d \tau d \H^2(x)\\
&=  \int_{\pa E_0} f(x) \big(u_{k+1} - u_{k}\big)\big(1 + H_{E_0}u_k + K_{E_0}u_k^2\big)  \, d \H^2(x) +o(h^2)\\
&= \int_{\pa E_0} f(x) J_{\pa E_0}\Phi_k  \frac{u_{k+1} - u_k}{|N_k|}  \, d \H^2(x) +o(h^2) , 
\end{split}
\]
where the last equality follows from the definition of $|N_k|$.  Above the integral $ \int_{u_k(x)}^{u_{k+1}(x)}$ is understood to be $-\int_{u_{k+1}(x)}^{u_k(x)}$ in the case $u_{k+1}(x) < u_{k}(x)$. 
We compute the same integral in a different way  by defining $\Psi_{k}(\cdot, \tau): \pa E_k \to \R^3$ as $\Psi_{k}(x, \tau) = x + \tau \nu_{E_k}(x)$, recall that $\Phi_k(x) = x + u_k(x)\nu_{E_0}$,
\[
\begin{split}
\int_{\R^3} &f(\pi_{\pa E_0}(x))(\chi_{E_{k+1}} - \chi_{E_{k}})\, dx\\
 &=  \int_{\pa E_k}  \int_{-r_0}^{r_0}  f(\pi_{\pa E_0}(\Psi_{k}(x, \tau)) \big(\chi_{E_{k+1}}(\Psi_{k}(x, \tau)) - \chi_{E_{k}}(\Psi_{k}(x, \tau)) \big) J_{\pa E_k}\Psi_{k}(\cdot, \tau) \, d \tau d \H^2(x) \\
&=   \int_{\pa E_k}  \int_{0}^{\psi_k(x)}  f(\pi_{\pa E_0}(\Psi_{k}(x, \tau))J_{\pa E_k}\Psi_{k}(\cdot, \tau) \, d \tau d \H^2(x)\\
&=   \int_{\pa E_k} \psi_k(x) f(\pi_{\pa E_0}(x))\, d \H^2(x) + o(h^2)   \\
&=  \int_{\pa E_0} \psi_k(\Phi_k(x))  f(x) J_{\pa E_0}\Phi_k  \, d \H^2(x) + o(h^2). 
\end{split}
\]
We compare the two formulas above and choose
\[
f(x)  =  \frac{1}{ J_{\pa E_0}\Phi_k} \left( \psi_k\circ \Phi_k  -  \frac{u_{k+1} - u_k}{|N_k|} \right)
\] 
and obtain  the estimate \eqref{eq:normal-velo-disc}.

Using the estimates from above we deduce that  we may write the approximate flat flow  $(E^{h,\delta}(t))_{t \geq 0}$ starting from  $E_0$ as 
\[
\pa E^{h,\delta}(t) = \{ x + u^{h,\delta}(x,t)\nu_{E_0}(x) : x \in \pa E_0 \}
\]
and the heightfunctions  $u^{h,\delta}$ satisfy
\[
\|u^{h,\delta}(\cdot ,t)\|_{H^4(\Sigma_0)} \leq C \quad \text{and} \quad \frac{1}{h}\| u^{h,\delta}(\cdot,t) - u^{h,\delta}(\cdot, t-h)  \|_{L^2(\Sigma_0)} \leq C
\]
for all $t \in [h,\sigma_0]$.  We may thus pass to the limit as $h \to 0$, up to a subsequence,  and deduce that the limit function $u(x,t) = \lim_{h \to 0} u^{h,\delta}(x,t)$ exists for every $t \in [0,\sigma_0]$. By \eqref{eq:normal-velo-disc}  the discrete normal velocity  converges, up to a subsequence,  to
\[
\lim_{h \to 0} v_k(\Phi_k(x))  = \frac{\partial_t u}{|N(x,t)|} \quad \text{in }\, L^2(E_0), \, t \in [0,\sigma_0],
\]
where $|N(x,t)|=  \frac{J_{\pa E_0}\Phi_t }{1 + H_{E_0}u(x,t) + K_{E_0}u^2(x,t)} $ and   $\Phi_t(x) = x + u(x,t) \nu_{E_0}$.  Let us show that the limit function solves the surface diffusion equation. 

To this aim we fix $\varphi \in C^2(\R^3)$,  multiply  the Euler-Lagrange equation  \eqref{eq:Euler-Lagrange-2}  with $\varphi$ and integrate by parts 
\[
\int_{\pa E_k^h}  \frac{\psi_k}{h}  \varphi  \,  d \H^2(x) = \int_{\pa E_k^h} \Delta_{\pa E_k^h} \varphi \,   ( H_{E_k^h} - \Delta_{\pa E_k^h} \psi_k + R_0)  \,   d \H^2(x) .
\]
Passing to the limit as $h \to 0$ and using the bounds  $\|\psi_k\|_{H^2} \leq h^{\gamma}$  and $\|R_0\|_{L^2} \leq h^{\gamma}$, which follow from the results in  the previous section,   yield 
\[
\begin{split}
\int_{\pa E_0}  &\frac{\partial_t u}{|N|}\varphi (\Phi_t(x))  J_{\pa E_0} (\Phi_t(x))   \,  d \H^2(x) = \int_{\pa E^{\delta}(t)} \Delta_{\pa E^{\delta}(t)} \, \varphi   H_{E^{\delta}(t)} \,   d \H^2(x)\\
&=  \int_{\pa E^{\delta}(t)}  \varphi  \, \Delta_{\pa E^{\delta}(t)} H_{E^{\delta}(t)} \,   d \H^2(x) =  \int_{\pa E_0}  \varphi (\Phi_t(x))  \, (\Delta_{\pa E^{\delta}(t)} H_{E^{\delta}(t)}) \,   J_{\pa E_0} (\Phi_t(x)) \,  d \H^2(x)  .
\end{split}
\]
Hence, we conclude that the limit flow $E^{\delta}(t)$, parametrized by the family of diffeomorphisms  $\Phi_t(x) = x + u(x,t) \nu_{E_0}$, is a strong solution to the surface diffusion equation. To be more precise, using the expansion of the mean curvature from Lemma \ref{lem:expand-mean-curv}
 and the expansion of the Laplace-Beltrami operator from e.g. \cite[formula (3.4)]{FJM3D} we deduce that $u : \Sigma_0 \times [0,\sigma_0] \to \R$, where $\Sigma_0 = \pa E_0$,  is a strong solution of the equation   (see \cite[formula (3.6)]{Antonia})
\begin{equation}\label{eq:surf-diff-pde}
\pa_t u = - \Delta_{\Sigma_0}^2 u + \la A(x,u, \bar \nabla u ), \bar \nabla^4 u \ra + J(x,u, \bar \nabla u, \bar \nabla^2 u, \bar \nabla^3 u ) \qquad \text{on }\, \Sigma_0 \times [0,\sigma_0],
\end{equation}
with initial datum $u(\cdot, 0) = 0$, where  the tensor field satisfies $A(x,0,0)= 0$ and $J$ is a smooth function with sublinear growth on $\bar \nabla^3 u $. By a strong solution we mean that $u \in W^{1,\infty}(0,\sigma_0; L^2(\Sigma_0))\cap L^\infty(0,\sigma_0; H^4(\Sigma_0))$ and it satisfies the equation \eqref{eq:surf-diff-pde} almost  everywhere. By  standard Gr\"onval  argument one may deduce that the strong solution to \eqref{eq:surf-diff-pde} with zero initial datum is unique. This proves that the limiting  flat flow coincides with the classical solution for the time interval $[0,\sigma_0]$ and the claim follows. 
\end{proof}

\begin{remark}
\label{rem:thm2}
We may quantify the statement of Theorem  \ref{thm:thm2} more precisely as follows.  Assume that  $E_0$ is $C^5$ regular, $E_0 \subset B_{R_0/2}$, satisfies the UBC with radius $2r_0$,
and the heightfunction in \eqref{eq:thm-2-height} satisfies 
\[
\|\psi_0 \|_{L^2(\pa E_0)} \leq L_0 h \quad \text{and} \quad \|\Delta_{\pa E_0}\psi_0 \|_{L^2(\Sigma_0)} \leq  L_0\sqrt{h}.
\]
Then there is  $K_0 = K_0(r_0,R_0, L_0)$ such that if
\[
\|\Delta_{\pa E_0} H_{E_0} \|_{L^2} \leq K_0  \quad \text{and} \quad \|\Delta_{\pa E_0} H_{E_0} \|_{H^1} \leq  K_0 h^{-\frac14}
\]
then  the approximate flat flow $(E^{(\delta,h)}(t))_{t\geq 0}$, where $\delta = \delta(r_0,R_0,K_0)$ is small, also satisfies  the UBC with radius $r_0$, $E^{(\delta,h)}(t) \subset B_{R_0}$ and 
\[
\|\Delta_{\pa E^{(\delta,h)}(t)} H_{E^{(\delta,h)}(t)} \|_{L^2} \leq  K_0  \quad \text{and} \quad \|\Delta_{\pa E^{(\delta,h)}(t)} H_{E^{(\delta,h)}(t)} \|_{H^1} \leq K_0 h^{-\frac14}
\]
for all $t \in [0,\sigma_0]$, where $\sigma_0 = \sigma_0(r_0,R_0,K_0)$.   
\end{remark}

\begin{remark}
\label{rem:uniqueness-strong}
We note that the strong solution of \eqref{eq:surf-diff-pde} with a zero initial datum $\psi(x,0) = 0$ is always unique. The argument in the proof of Theorem \ref{thm:thm2} implies that  if the approximate  flat flow $(E^{\delta, h}(t))_{t \geq 0}$ is such that $E^{\delta, h)}(0)$ is uniformly $C^5$-regular and $E^{\delta, h}(t)$  satisfies   the UBC with  radius $r_0$ and $\|\Delta_{\pa E^{\delta, h}(t)} H_{E^{\delta, h}(t)}\|_{L^2}\leq K_0$ for all $t \in [0,\sigma_0]$, then the limiting  flat flow coincides with the classical solution for the time interval $[0,\sigma_0]$.  

\end{remark}

\begin{proof}[\textbf{Proof of Theorem \ref{thm:thm1}}]
We assume that the classical solution   $(E_t)_{t \in [0,T_0)}$ exists for the time interval $[0, T_0)$ and fix  $T < T_0$. Since the classical solution is regular in $[0,T]$, there are  $r_0>0$ and $L_0>1$, depending not only on $E_0$ but also on $T$, such that $E_t$ satisfies the UBC with radius $4r_0$ and 
\begin{equation}
\label{eq:thm-1-0}
\|\Delta_{\pa  E_t} H_{ E_t}\|_{L^2}\leq  \frac{L_0}{2}\qquad \text{for all } t \in [0,T]. 
\end{equation} 
Also it holds  $ E_t \subset B_{R_0/2}$ for some $R_0=R_0(T)>0$ for all $t \leq T$. As in the proof of Theorem  \ref{thm:thm2}, the heightfunction $\psi_0$ as in \eqref{eq:thm-2-height} satisfies 
\[
\|\psi_0 \|_{L^2(\pa E_0)} \leq L_0 h \quad \text{and} \quad \|\Delta_{\pa E_0}\psi_0 \|_{L^2(\Sigma_0)} \leq  L_0\sqrt{h},
\]
for $L_0$ depending on $E_0$.  We choose    $K_0 = K_0(r_0,R_0,  L_0)$, $\sigma_0 = \sigma_0(r_0,R_0,K_0)$  and $\delta_0  = \sigma_0(r_0,R_0,K_0)$ as  in Remark~\ref{rem:thm2}, and denote $k_0 \in \mathbb{N}$ the number such that $\sigma_0 \in [k_0h, (k_0+1)h )$. Note that ultimately $K_0$, $\sigma_0$ and $\delta_0$ depend on $T$.

 Let $(E^{(\delta,h_i)}(t))_{t \geq 0}$ be the approximate flat flow associated with the chosen  flat flow $(E^\delta(t))_{t \geq 0}$, with $\delta <\delta_0$, and denote the associated sequence by $E_k^{(\delta,h_i)}$. We simplify the notation by $E^{h_i}(t) =  E^{(\delta,h_i)}(t)$ and $E_k = E_k^{(\delta,h_i)}$.  We claim that  for  every $t \in [0, T]$ it holds 
\begin{equation}
\label{eq:induction-P-l}
\begin{split}
E^{h_i}(t) \quad &\text{satisfies the UBC with radius }\, r_0, \quad E^{h_i}(t) \subset B_{R_0}\\
&\|\Delta_{\pa  E^{h_i}(t)} H_{E^{h_i}(t)}\|_{L^2}\leq K_0 \quad \text{and }\\
&\|\Delta_{\pa  E^{h_i}(t)} H_{E^{h_i}(t)}\|_{H^1}\leq K_0 h^{-\frac14},
\end{split}
\end{equation} 
for $h_i$ small. We point out that here we do not quantify the smallness of $h_i$ and we may have to pass to another subsequence of $h_i$.  Once we have proven \eqref{eq:induction-P-l}, the consistency  follows from Remark \ref{rem:uniqueness-strong}. 

Note that by Remark \ref{rem:thm2},  \eqref{eq:induction-P-l} holds for  all $t \in [0,\sigma_0]$. Let us choose any  $t_1 \in [\sigma_0,T)$ for which  \eqref{eq:induction-P-l} holds for all $t \leq t_1$. We will show that then \eqref{eq:induction-P-l} continues to hold for 
all $t \leq t_1+\frac{\sigma_0}{2}$. This will imply  that \eqref{eq:induction-P-l} holds for all $t \leq T$. 

Assume thus  that  \eqref{eq:induction-P-l} holds for  all $t \leq t_1$ and denote $k_1 $ the index for which $t_1- \frac{\sigma_0}{2} \in [k_1h_i, (k_1+1)h_i)$. Since $E_{k_1}$ satisfies the estimates  in \eqref{eq:induction-P-l}, we may use Proposition \ref{prop:apirori-estimates} and  Remark \ref{rem:thm2}  to deduce that  there is $\tilde \sigma = \tilde k h_i \geq \tilde c>0$, with $\tilde \sigma \leq \sigma_0$, such that for all $k \in [k_1,   k_1+ \tilde k]$ we may write the set $E_k$ as 
\[
\pa E_k  = \{ x + \psi_{k-1}(x) \nu_{E_{k-1}}(x) : x \in \pa E_{k-1} \}.
\]
We use  \eqref{eq:thm-2-1} to conclude 
\[
\| \psi_{k}\|_{L^2(\pa E_k)}^2 + h  \sum_{j=k_1}^{k_1+ \tilde k}  \| \Delta_{\pa E_j} \psi_{j}\|_{L^2(\pa E_j)}^2   \leq C h_i^2 \qquad \text{for all } \, k \in  [k_1,  k_1+ \tilde k], 
\]
for some constant $C$. Recall that $\tilde \sigma  =  \tilde k h_i \geq \tilde c >0$ for all $i$. Therefore we obtain from above that there is  an index $k_i  \in [k_1, k_1 +\tilde k/4]$ such that  
\begin{equation}
\label{eq:thm-1-1}
\| \psi_{k_i}\|_{L^2(\pa E_{k_i})} +  \| \Delta_{\pa E_{k_i} } \psi_{k_i} \|_{L^2(\pa E_{k_i})} \leq  C h_i.
\end{equation}

Next we recall that since the estimates in  \eqref{eq:induction-P-l}  hold for all $t \leq t_1$ and $k_i h_i \leq (k_1 +\tilde k/4)h_i \leq t_1$,  the set  $E_{k_i}$  satisfies the UBC with radius $r_0$ and is, in fact, uniformly $C^{2,\alpha}$-regular.
Let $\hat t_{i}  =  k_i h_i$ and define 
\[
v^{h_i}(x,\hat t_{i}) = \frac{\psi_{k_i}(x)}{h_i}. 
\]
Then \eqref{eq:thm-1-1} together with Lemma \ref{lem:hilbert-norm} and Sobolev embedding  imply that 
\[
 \| v^{h_i}(\cdot, \hat t_{i}) \|_{C^\alpha(\pa E^{h_i}(\hat t_{i}))} \leq  C .
\]
By choosing a subsequence of $h_i$ if needed, we may assume that 
\[
\lim_{h_i \to 0} \hat t_{i} =  \hat t \in [t_1- \frac{\sigma_0}{2} ,t_1].
\]
Therefore we deduce  by compactness  that 
\begin{equation}
\label{eq:thm-1-2}
\lim_{h_i \to 0 } \|v^{h_i}(\cdot, \hat t_{i})\|_{L^2(\pa E^{h_i}(\hat t_{i}))} = \|v(\cdot, \hat t)\|_{L^2(\pa E^\delta(t))},
\end{equation} 
by possible choosing another subsequence. Since we assume that \eqref{eq:induction-P-l} holds  for all $t \leq t_1$, and $\hat t \leq t_1$,  then by Remark  \ref{rem:uniqueness-strong} the flat flow agrees with the classical solution up to $t_1$. Then using \eqref{eq:normal-velo-disc} with $E_0$ replaced by $E_{t_1 - \tfrac{\sigma_0}{2}}$, we deduce that $v(\cdot, \hat t)$ agrees with the normal velocity $V_t$ of the classical solution $(E_t)$ at time $\hat t$, i.e.,  
\[
\|v(\cdot, \hat t)\|_{L^2(\pa E^\delta(\hat t))} = \|V_{\hat t}\|_{L^2(\pa  E_{\hat t})}. 
\]
Since $V_t = \Delta_{\pa E_t} H_{E_t}$ we deduce by \eqref{eq:thm-1-0} that 
\[
 \|V_{\hat t}\|_{L^2(\pa  E_{\hat t})} =  \|\Delta_{\pa E_{\hat t}} H_{E_{\hat t}}\|_{L^2(\pa  E_{\hat t})} \leq \frac{L_0}{2}.
\]
Therefore we conclude that it holds 
\[
\| \psi_{k_i}\|_{L^2(\pa E_{k_i})}  \leq  L_0 h_i,
\]
when $h_i$ is small. By  \eqref{eq:thm-1-1}  we also have
\[
\| \Delta_{\pa E_{k_i} } \psi_{k_i} \|_{L^2(\pa E_{k_i})} \leq  C h_i \leq L_0 \sqrt{h_i}. 
\]
when $h_i$ is small. Finally, since $E^\delta(\hat t) =  E_{\hat t}$ satisfies the UBC with radius $4 r_0$, then  by the uniform  $C^{2,\alpha}$-regularity we also deduce that $E_{k_i} =E^{h_i}(\hat t_i)$  satisfies the UBC with radius $2r_0$, when $h_i$ is small. 

We apply Remark \ref{rem:thm2} for $E_{k_i} =  E^{h_i}(\hat t_{i})$ in place of $E_0$ and deduce that the approximate flat flow satisfies
\[
\| \Delta_{\pa E^{h_i}(t)} H_{E^{h_i}(t)} \|_{L^2}\leq K_0\quad \text{and} \quad \|\Delta_{\pa  E^{h_i}(t)} H_{E^{h_i}(t)}\|_{H^1}\leq K_0 h^{-\frac14}
\]
for all $t \in [\hat t_i, \hat  t_i + \sigma_0]$.  Since $\hat t_i \in  [t_1- \frac{\sigma_0}{2} ,t_1]$,  this means that  the sets $E^{h_i}(t)$  satisfy  the bounds in \eqref{eq:induction-P-l} for all $t \leq t_1+\frac{\sigma_0}{2}$, when $h_i$ is small. 

We note that  by showing that the estimate \eqref{eq:induction-P-l} continues to hold from $t_1$ to $ t_1+\frac{\sigma_0}{2}$ we might have to pass to another subsequence from the original sequence $h_i$, but this is of minor importance as we need to repeat the   argument only finitely many times. Hence, the approximate flat flow $E^{h_i}(t)$ satisfies the bounds in  \eqref{eq:induction-P-l} for all $t \leq T$ and the claim  follows from Remark~\ref{rem:uniqueness-strong}. 
\end{proof}

\section*{Aknowledgments} 
During the completion of this paper V.J. and A.K were supported by the Academy of Finland grant 314227, while N.F. was supported by the Alexander von Humboldt Foundation.

\end{document}